%math.GT/0512594, 98ebw
%TOPOL-D-09-00094R1 http://ees.elsevier.com/topol/ skopenkov slovenija9
\input amstex
\documentstyle{amsppt}
\NoBlackBoxes
\magnification=\magstep1
%arxive
\advance\vsize0cm\voffset=-1.5cm\advance\hsize1.5cm\hoffset0cm
 \define\R{{\Bbb R}} \define\Z{{\Bbb Z}} \define\C{{\Bbb C}}
\def\pr{\mathop{\fam0 pr}}
\def\Sq{\mathop{\fam0 Sq}}

\def\id{\mathop{\fam0 id}}

\def\delet{\mathaccent"7017 }
\def\rel#1{\allowbreak\mkern8mu{\fam0rel}\,\,#1}
\def\Int{\mathop{\fam0 Int}}
\def\Con{\mathop{C}}
\def\Cl{\mathop{\fam0 Cl}}

\def\im{\mathop{\fam0 im}}

\def\t{\widetilde}

\def\tri{1.1\ }
\def\Tri{1.1}
\def\eff{1.2\ }
\def\Eff{1.2}
\def\is{1.3\ }
\def\Is{1.3}
\def\pri{1.4\ }
\def\Pri{1.4}

\def\Bh{2.1}

\def\Com{2.2}

\def\Rh{2.3}
\def\ru{2.4\ }
\def\Ru{2.4}

\def\Eu{2.5}

\def\Cor{2.6}

\def\Gso{3.1}

\def\Exn{3.2}
\def\bp{3.3\ }
\def\Bp{3.3}
\def\ret{3.4\ }
\def\Ret{3.4}

\def\Add{3.5}
\def\ntr{3.6\ }
\def\Ntr{3.6}
\def\adk{4.1\ }
\def\Adk{4.1}
\def\red{4.2\ }
\def\Red{4.2}
\def\fib{4.3\ }
\def\Fib{4.3}
\def\bor{4.4\ }
\def\Bor{4.4}
\def\clc{4.5\ }
\def\Clc{4.5}
\def\ale{4.6\ }
\def\Ale{4.6}

\newpage
\topmatter
\title A classification of smooth embeddings of 4-manifolds in 7-space, I \endtitle
\author Arkadiy Skopenkov \endauthor
\address
Independent University of Moscow, B. Vlasyevskiy, 11, 119002, Moscow, Russia, and
Faculty of Informatics and Calculation Technics, Moscow Institute of Physics and Technology, 141700, Dolgoprudnyi, Russia.
e-mail: skopenko\@mccme.ru \endaddress
\thanks
The author gratefully acknowledges the support from Deligne 2004 Balzan
prize in mathematics, the Russian Foundation for Basic Research Grants
07-01-00648a and 06-01-72551-NCNILa
\endthanks

\abstract
We work in the smooth category.
Let $N$ be a closed connected $n$-manifold and assume that $m>n+2$.
Denote by $E^m(N)$ the set of embeddings $N\to\R^m$ up to isotopy.
The group $E^m(S^n)$ acts on $E^m(N)$ by embedded connected summation of
a manifold and a sphere.
If $E^m(S^n)$ is non-zero (which often happens for $2m<3n+4$) then until
recently no complete readily calculable description of $E^m(N)$ or of this
action were known.
%(as far as I know).
Our main results are {\it examples of the triviality and the effectiveness of
this action}, and {\it a complete readily calculable isotopy classification of
embeddings into $\R^7$ for certain 4-manifolds $N$.}
The proofs use new approach based on the Kreck modified surgery theory
and the construction of a new invariant.

{\bf Corollary.}
{\it (a) There is a unique embedding $f:\C P^2\to\R^7$ up to isoposition
(i.e. for each two embeddings $f,f':\C P^2\to\R^7$ there is a diffeomorphism
$h:\R^7\to\R^7$
%(not necessarily orientation-preserving)
such that $f'=h\circ f$).

(b) For each embeddings $f:\C P^2\to\R^7$ and
%each non-trivial embedding
$g:S^4\to\R^7$ the embedding $f\#g$ is isotopic to $f$.}
\endabstract
\endtopmatter

\document

\head 1. Introduction \endhead

{\bf Knotting Problem for 4-manifolds.}

This paper is on the classical Knotting Problem: {\it given an $n$-manifold $N$
and a number $m$, describe isotopy classes of embeddings $N\to\R^m$}.
For recent surveys see [RS99, Sk08, HCEC].
We work in the smooth category.
%Let $ E^m(N)$ be the set of embeddings $N\to\R^m$ up to isotopy.
%(or, which is the same, up to orientation-preserving isoposition).

For $2m\ge3n+4$ there are some complete readily calculable classifications of
isotopy classes [Sk08, \S2, \S3, HCEC].
\footnote{By readily calculable classification I mean a classification
in terms of homology of a manifold (and certain structures in homology like
intersection, characteristic classes etc.).
A readily calculable classification is also a reduction to
calculation of stable homotopy groups of spheres (or to another standard
algebraic problem involving only homology of our manifold, which problem is
solved in particular cases, although could be unsolved generally).
An important feature of a useful classification is accessibility of statement
to general mathematical audience which is only familiar with basic notions of
the area.}
%Another important feature is whether the classification gives an algorithm
%and how complex the algorithm is, cf. [M?].
%in turn is the first approximation to beauty.
If
$$2m<3n+4$$
and a closed manifold $N$ is different from disjoint union of homology spheres,
then until recently no complete readily calculable descriptions of isotopy
classes were known (up to my knowledge),
in spite of the existence of interesting approaches of Browder-Wall and
Goodwillie-Weiss [Wa70, GW99, CRS04].
\footnote{The approach of [GW99] gives a modern abstract proof of certain
earlier known results; I am grateful to M. Weiss for indicating that the
approach also gives explicit results on higher homotopy groups of the space of
embeddings $S^1\to\R^n$.}
For recent results see [Sk06] (a classification of embeddings
$S^p\times S^q\to\R^m$), [Sk08'] (a classification of embeddings of 3-manifolds
into $\R^6$), [CS], [CS11] (see comments on this closely related paper after Lemma
\Is),
[CRS07, CRS08] (rational classification).
For {\it piecewise linear} classification see [Sk07, Sk08, \S2, \S3].

In particular, a complete readily calculable classification of embeddings of a
closed connected 4-manifold $N$ into $\R^m$ was only known either for $m\ge8$
or for $N=S^4$ and $m=7$ by Wu, Haefliger, Hirsch and Bausum:
%(or for certain 4-manifolds and $m=7$ in the piecewise-linear category by
%Bo\'echat-Haefliger and the author):
$$\# E^m(N)=1\quad\text{for}\quad m\ge9,$$
$$E^8(N)=\cases H_1(N;\Z_2) & N\text{ orientable}\\
\Z\oplus\Z_2^{s-1} & N\text{ non-orientable and }H_1(N;\Z_2)\cong\Z_2^s
\endcases,$$
$$E^7(S^4)\cong\Z_{12}.$$
%$$E_{PL}^7(N)=H_2(N;\Z)\quad\text{if}\quad H_1(N;\Z)=0.$$
%$$E^7_{PL}(S^1\times S^3)\cong\Z\oplus\Z\oplus\Z_2.$$
Here $E^m(N)$ is the set of smooth embeddings $N\to\R^m$ up to smooth isotopy;
the equality sign between sets denotes the existence of a bijection; the
isomorphism is a group isomorphism for the `connected sum' group structure
on $E^7(S^4)$ [Ha66].
See references in [Sk08, \S2, \S3, HCEC]; cf. [CS11, CS].
\footnote{The known {\it existence} results for closed 4-manifolds $N$ are as
follows:
\newline
$\bullet$ $N$ embeds into $\R^8$;
\newline
$\bullet$ if $N$ is orientable, then $N$ embeds into $\R^7$
[Hi65, Fu94, cf. BH70, Fu02];
\newline
$\bullet$ $N$ embeds into $\R^7$ if and only if $\overline W_3(N)=0$ [Fu94];
\newline
$\bullet$ an orientable $N$ PL embeds into $\R^6$ if and only if $\overline w_2(N)=0$
[Ma80, Corollary 10.11, CS79].
\newline
$\bullet$ an orientable $N$ smoothly embeds into $\R^6$ if and only if
$\overline w_2(N)=0$ and $\sigma(N)=0$ [Ma80, Corollary 10.11, CS79, Ru82].}

%Cf. [GoSt, \S].CHR01
%Concerning embeddability of 4-manifolds into $\R^5$ see [Coh84I, Coh84T, CHS02].)
%(DIFF or PL: recall that any PL structure on a 4-manifold admits a unique smoothing).

One of the main results of this paper (the Triviality Theorem \Tri(b)) is
{\it a complete readily calculable classification of embeddings of certain
4-manifolds into $\R^7$}, cf. [CS11, CS].
%A non-trivial corollary can be stated without any notation.

\bigskip
{\bf Main results.}

Consider the `connected sum' group structure on $E^m(S^n)$ [Ha66].
By [Ha61, Ha66, Corollary 6.6, Sk08, \S3],
$E^m(S^n)=0\quad\text{for}\quad 2m\ge3n+4.$
However, $E^m(S^n)\ne0$ for many $m,n$ such that $2m<3n+4$,
\footnote{This differs from the Zeeman-Stal\-lings Unknotting Theorem:
{\it for $m\ge n+3$ any PL or TOP embedding $S^n\to S^m$ is PL or TOP
isotopic to the standard embedding}.}
e.g. $E^7(S^4)\cong\Z_{12}$.
\footnote{This follows from [Ha66, 4.11, cf. Ha86] and the well-known Lemma \Gso.}

In this and the next subsections we assume that $N$ is a closed $n$-manifold and
$m\ge n+3$.
The group $E^m(S^n)$ acts on the set $E^m(N)$ by connected summation of
embeddings $g:S^n\to\R^m$ and $f:N\to\R^m$ whose images are contained in
disjoint balls.
\footnote{Since $m\ge n+3$, the connected sum is well defined, i.e. does not
depend on the choice of an arc between $gS^n$ and $fN$.
If $N$ is not connected, we assume that a component of $N$ is chosen and
we consider embedded connected summation with this chosen component. }
Various authors have studied analogous connected sum action of the group of
homotopy $n$-spheres on the set of $n$-manifolds topologically homeomorphic to
given manifold [Le70].

The quotient of $E^m(N)$ modulo the above action of $E^m(S^n)$ is known
in some cases.
\footnote{In those cases when this quotient coincides with the set of PL
embeddings $N\to\R^m$ up to PL isotopy
and when the latter set was known [Sk08, \S2, Sk02, Sk07].}
Thus in these cases the Knotting Problem is reduced to the description of
the above action of $E^m(S^n)$ on $E^m(N)$.
Until recently {\it no results were known} on this action for $E^m(S^n)\ne0$
and $N$ not a disjoint union of homology spheres.
For recent results see [Sk06, Smoothing Theorem, Sk08', CS11, CS];
for $m=n+2$ see [Vi73].

The main results of this paper are the following {\it examples of the
triviality and the effectiveness of the above action for embeddings of
4-manifolds into $\R^7$}.

We omit $\Z$-coefficients from the notation of (co)ho\-mo\-lo\-gy groups.

\smallskip
{\bf The Triviality Theorem \Tri.}
{\it Let $N$ be a closed connected smooth 4-manifold
such that $H_1(N)=0$ and the signature $\sigma(N)$ of $N$ is free of squares
(i.e. is not divisible by the square of an integer $s\ge2$).

(a) For each embeddings $f:N\to\R^7$ and $g:S^4\to\R^7$ the embedding $f\#g$ is
isotopic to $f$ (although $g$ could be non-isotopic to the standard embedding).

(b) There is a 1--1 correspondence}
$$\varkappa : E\phantom{}^7(N)\to
\{x\in H_2(N)\ |\ x\mod2=PDw_2(N),\ x\cap x=\sigma(N)\}.$$

Here $PDw_2(N)$ is Poincar\'e dual of the 2nd Sitefel-Whitney class, and $\cap$
is the intersection product, cf. Remark \Rh.

E. g. $N=\C P^2$ satisfies the assumption of the Triviality Theorem \Tri,
so $ E^7(\C P^2)$ is in 1--1 correspondence with
$\{+1,-1\}\subset\Z\cong H_2(\C P^2)$.
\footnote{The two isotopy classes of embeddings $\C P^2\to\R^7$ are represented
by the standard embedding and by its composition with the reflection of $\R^7$.
The standard embedding is given by
%the well-known formula
$(x:y:z)\mapsto$
\linebreak
$\mapsto(x\overline y, y\overline z, z\overline x,2|x|^2+|y|^2)$, where $|x|^2+|y|^2+|z|^2=1.$
Alternatively, see [BH70, p. 164, Sk08', \S5].}
% See also [Ca39, KK95, Ma75].

%(This standard embedding is a part of a linking $\C P^2\sqcup \C P^2\to S^7$
%such that each component is a deformation retract
%of the other, cf. the Complement Lemma \com below.) }

\smallskip
{\bf The Effectiveness Theorem \Eff.}
{\it Let $N$ be a closed smooth simply-connected 4-manifold
%such that $\Sigma N$ retracts to $\Sigma N_0$.
and $f_0:N\to\R^7$ an embedding such that $f_0N\subset\R^6$.
%(thus $N$ is embeddable into $\R^6$).
Then each for pair of non-isotopic embeddings $g_1,g_2:S^4\to\R^7$
%the
embedding $f_0\#g_1$ is not isotopic to $f_0\#g_2$.}

\smallskip
{\bf Remark.}
Let $N$ be a closed connected orientable smooth 4-manifold.
Then the following conditions are equivalent:

(1) $N$ is embeddable into $\R^6$;

(2) $w_2(N)=0$ and $\sigma(N)=0$;

(3) the normal bundle of each embedding $f:N\to S^7$ is trivial.

For simply-connected $N$ each of these conditions is equivalent to

(4) the intersection form of $N$ is that of $\#_i(S^2\times S^2)$;
\footnote{The connected sum in (4), (5) and (6) could have zero summands; then the connected sum is $S^4$.}

(5) $N$ is homotopy equivalent to $\#_i(S^2\times S^2)$;

(6) $N$ is topologically homeomorphic to $\#_i(S^2\times S^2)$.

%$N$ is parallelizable No:e=0;
%The latter implies that $\Sigma N$ retracts to $\Sigma N_0$.}

\smallskip
{\it Proof.} (1) is equivalent to (2) by [GS99, Theorem 9.1.21 and Remark
9.1.22, cf. CS79, Ru82].
(2) is equivalent to (3) by the Dold-Whitney Theorem [DW59],
cf. [CS11, the Normal Bundle Lemma].

By the Whitehead and the Freedman Theorems, (4) is equivalent to (5) and (6).

Clearly, (4) implies (2).
For simply-connected $N$ (2) implies that the intersection form of $N$ is
indefinite and even, so (4) holds by [Ma80, Theorem 1.9.2].
\qed

\bigskip
{\bf Ideas of proof.}
\nopagebreak

The Effectiveness Theorem \eff is proved in \S3 (cf. Theorem \Ntr).
The proof is based on the construction of a new {\it attaching invariant} for
certain embeddings $N\to\R^m$, generalizing the Haefliger-Levine attaching
invariant of embeddings $S^n\to\R^m$.

The proof of the Triviality Theorem \tri is much more non-trivial.
The proof is based on the following idea,
which is useful not only in these dimensions [Sk08'] and not only to
describe the action of $E^m(S^n)$ on $E^m(N)$ [FKV87, FKV88].
%This idea was already used to construct exotic knottings of 2-surfaces in $S^4$

Fix orientations on $N$ and on $\R^m$.
For an embedding $f:N\to\R^m$ denote by

$\bullet$ $C_f$ the closure of the complement in $S^m\supset\R^m$ to a tubular
neighborhood of $fN$,

$\bullet$ $\nu_f:\partial C_f\to N$ the restriction of the normal bundle of $f$.

%An orientation-preserving diffeomorphism
%$\varphi:\partial C_f\to\partial C_{f'}$ such that $\nu_f=\nu_{f'}\varphi$
%is called simply an {\it isomorphism}.

\smallskip
{\bf Lemma \Is.}
{\it For a closed connected manifold $N$ embeddings $f,f':N\to\R^m$ are
isotopic if and only if there is an orientation-preserving bundle isomorphism
$\varphi:\partial C_f\to\partial C_{f'}$ which extends to an
orientation-preserving diffeomorphism $C_f\to C_{f'}$.}

\smallskip
{\it Proof.}
The `only if' part is obvious, so let us prove the `if' part.
The isomorphism $\varphi$ also extends to an orientation-preserving
diffeomorphism $S^m-\Int C_f\to S^m-\Int C_{f'}$.
%Therefore $\Sigma\cong S^m\#\Sigma\cong S^m$. $C_f\cong C_{f'}$ and
Hence $\varphi$ extends to an orientation-preserving diffeomorphism
$\R^m\cong \R^m$.
Since any orientation-preserving diffeomorphism of $\R^m$ is isotopic to
the identity, $f$ and $f'$ are isotopic.
\qed

\smallskip
So results on the Diffeomorphism Problem can be applied to Knotting Problem.
In this way there were obtained embedding theorems in terms of Poincar\'e
embeddings [Wa70].
But 'these theorems reduce geometric problems to algebraic problems which are
even harder to solve' [Wa70].
One of the main problems is that in general (i.e. not in simpler cases like
that of the Effectiveness Theorem \Eff) it is hard to work with the homotopy
type of $C_f$ (which is sometimes unknown even when the classification of
embeddings is known).
The main idea of our proof is to apply the modification of surgery [Kr99]
which allows to classify $m$-manifolds using their homotopy type just below
dimension $m/2$.

The relation of this paper to a closely related paper [CS11] is as follows.
Shortly, the proof for the simplest (but non-trivial) cases (like $\C P^2$ in
$\R^7$) is presented in this paper, while more complicated proof for more
general case is given in [CS11].
More precisely, the main result of [CS11] is a description of $ E^7(N)$ for
each closed connected 4-manifold $N$ such that $H_1(N)=0$.
This result recovers the Triviality Theorem \tri completely and the
Effectiveness Theorem \eff only for the case $H_1(N)=0$.
The proof in this paper is almost disjoint from the proof in [CS11] and is much shorter.
The difference in applying the modification of surgery [Kr99]
is that here we use $BO\left<5\right>\times\C P^\infty$-surgery while in [CS11]
$BSpin\times\C P^\infty$-surgery is used; the attaching invariant is not used
in [CS11].

In the first subsection of \S3 we present the well-known definition of the
attaching invariant $a: E^7(S^4)\to\Z_{12}$.
This subsection is formally not used later in \S3, where we generalize this
definition.
In \S4 we give a new proof of its injectivity based on [Kr99].
%This material uses only the definition of $C_f$ from the rest of the
%Introduction.
Except for this relation, \S3 and \S4 are independent on each other.

The Triviality Theorem \Tri(a) follows by the Bo\'echat-Haefliger Theorem
\Bh(a) and the Complement Lemma \Com(b) below, together with the following
result, which is the new and the most important part of the proof.

\smallskip
{\bf The Primitivity Theorem \Pri.}
{\it Let $N$ be a closed connected smooth 4-manifold and $f:N\to\R^7$ an
embedding such that $\pi_3(C_f)=0$ (and hence $H_1(N)=0$).
Then for each embedding $g:S^4\to\R^7$ the embedding $f\#g$ is isotopic to
$f$.}

\bigskip
%\smallskip
{\bf Acknowledgements.}
These results are based on ideas of and discussions with Matthias Kreck.
%(he could not be persuaded to be a coauthor of this paper)
They were originally announced in [KS05]
and presented at the 4th European
Congress of Mathematicians, at the Conference 'Topology, Geometry and
Combinatorics' (Stockholm, 2004), at the L. V. Keldysh Centennial
Conference (Moscow, 2004) and at the International Congress on
Differential Geometry (M\"unster, 2006).
I would like to acknowledge D. Crowley, P. Lambrechts, S. Melikhov, D. Tonkonog and
anonymous referees for useful discussions.

\head 2. Preliminaries  \endhead

%\bigskip
{\bf The Bo\'echat-Haefliger invariant.}
\nopagebreak

From now on, unless otherwise stated, we assume that

{\it $N$ is a closed connected orientable 4-manifold and $f:N\to\R^7$ is an
embedding}.

The {\it homology Alexander duality isomorphism} $A_f:H_i(N)\to H_{i+1}(C_f,\partial C_f)$ is
the inverse to the composition $H_5(C_f,\partial C_f)\to H_4(\partial C_f)\to H_4(N)$ of
the boundary map and the normal bundle map.
(This composition equals to the composition
$H_4(N)\to H^2(C_f)\to H_5(C_f,\partial C_f)$ of the Alexander
and Poincar\'e-Lefschetz duality isomorphisms [Sk08', the Alexander Duality
Lemma].)

The {\it homology Seifert surface} for $f$ is $A_f[N]$.

Define {\it the Bo\'echat-Haefliger invariant}
$$\varkappa :E^7(N)\to H_2(N)\quad\text{by}\quad \varkappa (f):=A_f^{-1}(A_f[N]\cap A_f[N]).$$
(A definition of $\cap$ is recalled in Remark \Rh.)

This new definition is equivalent to the original one [BH70] by the Section Lemma \Eu.

\smallskip
{\bf The Bo\'echat-Haefliger Theorem \Bh.}
%{\it Let $N$ be a closed connected orientable 4-manifold.
$$(a)\qquad
\im \varkappa \ =\ \{x\in H_2(N)\ |\ x\mod2=PDw_2(N),\ x\cap x=\sigma(N)\}.$$
\quad
(b) {\it If $H_1(N)=0$, then two embeddings $N\to\R^7$ with the same
$\varkappa$-invariant differ by a connected sum with an embedding $S^4\to\R^7$.}

\smallskip
E. g. $H_2(S^2\times S^2)\cong\Z\oplus\Z\supset\{(2k,2l)\ |\ kl=0\}=\im \varkappa $
for $N=S^2\times S^2$.
%examples of embeddings, and also that of CP^2

Here part (a) follows by the Section Lemma \Eu(b) and [BH70, Theorem 2.1], cf.
[Fu94].
Part (b) follows by the Section Lemma \Eu(b),
[BH70, Theorem 1.6] and smoothing theory [BH70, p. 156], cf. [Ha67, Ha].

%Note that the assumption $H_1(N)=0$ is essential in part (b) by [CS11].

The Triviality Theorem \Tri(b) follows from the Triviality Theorem \Tri(a) and
the Bo\'echat-Haefliger Theorem \Bh.

\smallskip
{\bf The Complement Lemma \Com.}
{\it If $H_1(N)=0$, then

(a) $C_f\simeq C_{\varkappa (f)}:=S^2\cup_{\varkappa (f)}(D^4_1\sqcup\dots\sqcup D^4_{b_2(N)})$.
Here we identify $H_2(N)$ and $\Z^{b_2(N)}$ by any isomorphism, so
$\varkappa (f)$ is identified with an ordered set of $b_2(N)$ integers, which set
defines a homotopy class of maps
$\partial(D^4_1\sqcup\dots\sqcup D^4_{b_2(N)})\to S^2$.
\footnote{This ordered set depends on the identification of $H_2(N)$ and
$\Z^{b_2(N)}$, but the homotopy type of $C_{\varkappa (f)}$ does not.
The homotopy equivalence $C_f\simeq C_{\varkappa (f)}$ is not canonical.}

%(so $C_{\varkappa (f)}$ is not a useful model for a complement, cf. proof of the
%Primitivity Theorem \Pri).

(b) $\pi_3(C_f)\cong\Z/d\Z$, where $d=0$ for $\varkappa (f)=0$ and}
$$d:=\max\{k\in\Z\ |\ \text{there is }y\in H_2(N):\ \varkappa (f)=ky\}
\quad\text{for}\quad \varkappa (f)\ne0.$$
%is the divisibility of $x$, i.e. $d(0):=0$ and
%$$d(x):=\max\{d\in\Z\ |\ \text{there is }x_1\in H_2(N):\ x=dx_1\}
%\quad\text{for}\quad x\ne0.$$}

%\smallskip
{\it Proof.} Part (b) follows by (a).
Let us prove (a).
By general position $C_f$ is simply-connected.
Since $H_1(N)=0$, by Alexander duality
$$H^2(C_f)\cong\Z,\quad H^4(C_f)\cong H_2(N)\quad\text{and}\quad H^i(C_f)=0
\quad\text{for}\quad i\ne2,4.$$
Now the Complement Lemma \Com(a) is obtained by taking $p=1$ in the following
statement implied by `homology decomposition of Eckmann-Hilton' [EH59].

{\it If $Y$ is a finite simply-connected cell-complex such that $H^i(Y)=0$ for
$i\ne0,2,4$, $H^2(Y)\cong\Z^p$, $H^4(Y)\cong\Z^q$ and $A$ is the
$(p\times q)$-matrix of the cup square $H^2(Y)\to H^4(Y)$ in some bases of
$H^2(Y)$ and $H^4(Y)$, then
$Y \simeq (S^2_1\vee\dots\vee S^2_p) \cup_A (e^4_1\sqcup\dots\sqcup e^4_q)$.}
\qed
\footnote{
{\it An alternative direct proof of the Complement Lemma \Com(b).}
Take the map $h_f:C_f\to\C P^\infty$ defined in \S4 at the beginning of the
proof of the Primitivity Theorem \Pri.
Then
$$\pi_3(C_f)\cong\pi_4(\C P^\infty,C_f)\cong H_4(\C P^\infty,C_f)
\cong H_4(\C P^\infty)/h_{f,*}H_4(C_f)\cong\Z/d\Z.$$
Here the fourth equality follows because for the dual map
$h_f^*:H^4(\C P^\infty)\to H^4(C_f)$ and the generator $a\in H^2(\C P^\infty)$
we have $h_f^*(a\cup a)=h_f^*a\cup h_f^*a=PDA_f[N]\cup PDA_f[N]$.
\qed
\newline
The Primitivity Theorem \pri and the Complement Lemma \Com(b) imply the
following.
\newline
$\bullet$ If $H_1(N)=0$ and $y\in H_2(N)$ is primitive (i.e. there are no
integers $d\ge2$ and elements $x\in H_2(N)$ such that $y=dx$), then
$\#\varkappa ^{-1}y$ is 0 or 1.
\newline
This could be proved analogously to the Primitivity Theorem
but using $\varkappa (f)=\varkappa (f')$ and [CS11, Agreement Lemma] instead of $f'=f\#g$.
This corollary and the Bo\'echat-Haefliger Theorem \Bh(a) give a proof of the
Triviality Theorem \tri without reference to the Bo\'echat-Haefliger Theorem
\Bh(b).
\newline
$\bullet$ Under the assumptions of the Primitivity Theorem \pri
the number of isotopy classes of smooth embeddings $f:N\to S^7$ for which
$\pi_3(C_f)\cong0$ equals to the number of primitive elements in $\im \varkappa $.}

%$\bullet$ the first one by the homotopy exact sequence of $(\Map h,C)$ because
%$\Map h\simeq\C P^\infty$,
%$\bullet$ the second one by the Hurewicz theorem,
%$\bullet$ the third one by the homology exact sequence of pair $(\C P^\infty,C)$ and
%$\bullet$ the fourth one because

\bigskip
{\bf Section Lemma \Eu.}

In this subsection we omit index $f$ from the notation.

\smallskip
{\bf Remark \Rh.} In this paper we mostly use the language of homology rather
than cohomology.
This makes the arguments more visual and so is (within geometric problems like
those treated here) more convenient to understand, check and apply the results.
Let us recall the main definitions (they are well known, see [Fe83], pp. 24-30,
where also some more details are given).
We give equivalent definitions in cohomological language
%(for the sake of the reader more used to cohomology).
(because all the required properties, if not found in the literature, can
either be derived from the cohomology properties or proved directly analogously
to them).
Cf. [MC06].

Let $Q$ be a compact smooth $q$-manifold.
Let $T$ be a smooth cell-decomposition of $Q$ in the sense of [RS72].
Denote by $H_i(T)$ the corresponding cellular homology groups.
Recall that $H_i(Q):=H_i(T)$ is independent of $T$.
Analogously one defines $H_i(Q,\partial Q)$ which we shortly denote by
$H_i(Q,\partial)$.

%Recall that $PDw_2(Q)\in H_{q-2}(T,\partial T)\cong H_{q-2}(Q,\partial)$
%is the obstruction to the existence of a family of $q-1$ linearly
%independent tangent to $Q$ vector fields over 2-skeleton of $T^*$.
%(In fact, this is used only in the statements of the Triviality and the
%Bo\'echat-Haefliger Theorems \tri and \Bh, not in the proofs.)

Let $T^*$ be the dual cell-decomposition.
Let $\overline T$ be the barycentric subdivision of $T$.
%then $\overline T$ is a refinement of $T^*$.
The intersection product $H_i(T)\times H_j(T^*)\to H_{i+j-q}(\overline T)$ is
defined using chain intersections.
This gives the intersection product $\cap:H_i(Q)\times H_j(Q)\to H_{i+j-q}(Q)$.
Analogously one defines the intersection products
$\cap_\partial:H_i(Q,\partial)\times H_j(Q)\to H_{i+j-q}(Q)$ and
$\cap_{\partial\partial}:
H_i(Q,\partial)\times H_j(Q,\partial)\to H_{i+j-q}(Q,\partial)$.

Denote Poincar\'e-Lefschetz duality (in any $q$-manifold $Q$) by
$$PD:H^i(Q)\to H_{q-i}(Q,\partial)\quad\text{and}
\quad PD:H^i(Q,\partial)\to H_{q-i}(Q).$$
Product $\cap':H^i(Q)\times H_j(Q)\to H_{j-i}(Q)$ is
%$\cup:H^i(Q)\times H^j(Q)\to H^{i+j}(Q)$ are
defined e.g. in [Pr07, 8.1].
Analogously one defines
$\cap'_\partial:H^i(Q,\partial)\times H_j(Q)\to H_{j-i}(Q)$ and
$\cap'_{\partial\partial}:H^i(Q)\times H_j(Q,\partial)\to H_{j-i}(Q,\partial)$.
%$\cap'_{\partial\partial}:H^i(Q,\partial)\times H_j(Q,\partial)\to H_{j-i}(Q)$
Clearly, $x\cap y=PDx\cap'_\partial y$,
$x\cap_{\partial\partial}y=PDx\cup_{\partial\partial}y$ and
$x\cap_\partial y=PDx\cap'y$.
\footnote{The intersection product $\cap$ is well defined, i.e. independent of
the triangulation, because $x\cap y=PDx\cap'_\partial y$ and $\cap'_\partial$
is well defined. (For details of this proof one uses the `chain-level
Poincar\'e isomorphism' $PD:C_i(T^*)\to C^{q-i}(T)$.)
Analogously by passing to cohomology one proves that the other intersection
products, the preimage homomorphism $s^!$ and the homology Euler class $PDe(p)$
below are well defined. }

In the sequel all the products $\cap,\cap_\partial,\cap_{\partial\partial}$ are
denoted simply by $\cap$, the domain of $\cap$ being clear from the context.

%(Passing to cohomology is the usual way to prove the independence of homology
%intersection products on the choice of triangulation.)

Let $p:E\to Q$ be the $D^k$-bundle associated to a real oriented
$k$-dimensional vector bundle.
Let $s_*$ be the zero section.
The {\it `preimage' homomorphism} $s^!:H_i(E)\to H_{i-k}(Q)$ is defined as
follows.
Take smooth cell-decompositions $T_B,T_E$ of $B,E$ such that $s$ is cellular.
Represent a class $x\in H_i(E)$ as a cellular cycle in the dual cell
decomposition to $T_E$.
Define $s^![x]$ as the $s$-preimage of $x$.
Clearly, $s^!=PD\circ s^*\circ PD^{-1}$.
{\it The homology Euler class} of $p$  is defined as
$PDe(p):=p_*(s_*[Q]\cap s_*[Q])=s^!s_*[Q]\in H_{q-k}(Q,\partial)$.
\footnote{This is clearly equivalent to one of the cohomological definitions
$e(p):=s^*s_!PD[Q]\in H^k(E)$.
The equality $p_*(s_*[Q]\cap s_*[Q])=s^!s_*[Q]$ is well known; here is a proof.
Represent $s_*[Q]$ as a cellular cycle $q$ in the dual cell decomposition to
$T_E$.
Identify $Q$ with $s(Q)$ by the embedding $s$.
Then both $p_*(s_*[Q]\cap s_*[Q])$ and $s^!s_*[Q]$ are represented by
the chain intersection of $q$ with the fundamental chain of $s(Q)$.
\newline
Alternative proofs could be obtained either passing to the smooth category and
using Thom transversality theorem, or using cohomological definition.}
%, or using [BRS76], pp. 25-26. Cf. [EG07], Theorem 1.1, [Me09], Lemma 2.1.}
%\footnote{It was suggested by a referee that the following proof of
%$p_*(s_*[Q]\cap s_*[Q])=s^!s_*[Q]$ is obtained using [BRS76], chapter II.
%By definition in p. 26, $e(p)=s^*s_*[Q]$.
%On the other hand, $s_*[Q]\cup s_*[Q]$ is $s_!s^*s_*[Q]$ by definition in
%p. 25.
%Hence $s_*[Q]\cup s_*[Q]=s_!e(p)$.
%Applying $p_!$ to both sides and observing that $ps=\id$,
%we obtain $p_!(s_*[Q]\cup s_*[Q])=e(p)$ as required.}

%$PDe(p)\in H_{q-k}(T,\partial T)\cong H_{q-k}(Q,\partial Q)$ being the
%obstruction to the existence of non-zero section over $k$-skeleton of $T^*$.)

\smallskip
{\it Definition of a weakly unlinked section.}
Let $N_0:=\Cl(N-B^4)$, where $B^4$ is a closed 4-ball in $N$.
Let $\zeta:N_0\to \nu^{-1}N_0$ be a section of the normal bundle
$\nu^{-1}N_0\to N_0$.
(This exists because $e(\nu)=0$.)
Consider the following diagram.
$$\minCDarrowwidth{0pt}\CD
H_4(N_0,\partial) @>> \zeta_* > H_4(\nu^{-1}N_0,\partial) @<< e <
H_4(\partial C,\nu^{-1}B^4) @<< j <  H_4(\partial C) @>> i > H_4(C)\endCD.$$
Here $j$ is the isomorphism from the exact sequence of pair, $e$ is the
excision isomorphism and $i$ is induced by the  inclusion.
Section $\zeta$ is called {\it weakly unlinked} if $ij^{-1}e^{-1}\zeta_*=0$.

\smallskip
{\bf Remark \Ru.}
In the definition of a weakly unlinked section we can replace $i$ by
$i':H_4(\partial C)\to H_4(S^7-fN_0)$.
\footnote{So the definition is equivalent to the following original definition
[BH70].
Denote by $|\cdot,\cdot|$ a distance function in $N$ such that $B^4$ is a ball
of radius 2.
Define a map
$$\overline\zeta:N\to S^7-fN_0\quad\text{by}\quad
\overline\zeta(x)=\cases \zeta(x)&x\in N_0\\ f(x)&|x,N_0|\ge1\\
|x,N_0|f(x)+(1-|x,N_0|)\zeta(x)&|x,N_0|\le1\endcases.$$
Section $\zeta$ is called {\it weakly unlinked} if
$\overline\zeta_*[N]=0\in H_4(S^7-fN_0)$.}
%or, equivalently, $A_{f,\xi}=0$
Indeed, let $\widehat\nu:S^7-\Int C\to N$ be the disk normal bundle.
The remark follows because inclusions induce isomorphisms
$H_4(C)\to H_4(C\cup\widehat\nu^{-1}B^4)\to H_4(S^7-fN_0)$.
(The first inclusion is an isomorphism by the exact sequence of pair,
the second because it is an inverse to a strong deformation retraction.)

\smallskip
{\bf Section Lemma \Eu.} {\it If $\zeta$ is a weakly
unlinked section, then

(a) $ej\partial A[N]=\zeta_*[N_0]$.

(b) $\varkappa (f)=PDe(\zeta^\perp)=\zeta^!ej\partial A[N]$, where $\zeta^\perp$ is the
oriented $S^1$-bundle that is the orthogonal complement to $\zeta$ in
$\nu|_{N_0}$, and for $k\ne0$ we identify $H_k(N)$ with $H_k(N_0,\partial)$
by the composition
$H_k(N)\overset{j_N}\to\to H_k(N,B^4)\overset{e_N}\to\to H_k(N_0,\partial)$
of the isomorphism from the exact sequence of pair and the excision
isomorphism.}

\smallskip
{\it Proof.} First we prove (a).
\footnote{Cf. [Sk08'], Unlinked Section Lemma (c); in our case a weakly
unlinked section need not extend to a section over $N$.
{\it An alternative proof of (a) using the original definition [BH70]} is as follows.
Take a smooth triangulation of $S^7$ such that $fN_0$, $fB^4$, $\partial C$,
$\zeta N_0$ and the union $b$ of segments $f(x)\overline\zeta(x)$, $x\in N$,
are subcomplexes.
%Denote by $\widehat\nu^{-1}N_0$ the tubular neighborhood of $fN_0$ in $S^7$
%(so that $S^7-\Int C-\widehat\nu^{-1}N_0\cong \Int B^4\times D^3$).
Since $\zeta$ is weakly unlinked and
$S^7-fN_0\simeq C\cup\widehat\nu^{-1}B^4$, there is a 5-chain $a$ in
$C\cup\widehat\nu^{-1}B^4$ such that $\partial a$ is represented by
$\overline\zeta N$.
(This 5-chain $a$ is in some refinement of the above triangulation, which
refinement is used in the rest of this proof.)
Recall that 5-chain $a\cap C$ is defined as the sum of simplices of $a$ that
are in $C$.
Recall that 5-chain $a\cap(S^7-C)$ is defined as the sum of closures of
simplices of $a$ whose interiors are in $S^7-C$.
(These definitions are of course different from the definition of the
intersection in homology.)
The sum of $a\cap(S^7-C)$ and the 5-chain in $S^7-\Int C$ represented by
$b$ is a homology between $\partial(a\cap C)$ and the 4-chain in
$S^7-\Int C$ represented by $fN$ in $S^7-\Int C$.
Then $\nu_*\partial[a\cap C]=[N]$, so by homology Alexander duality
$A=[a\cap C]$.
Hence $ej\partial A[N]=[(\nu^{-1}N_0)\cap\partial(a\cap C)]=\zeta_*[N_0]$
(here the intersection of a subset and a chain is defined analogously to
the above).}
Since $\zeta$ is weakly unlinked, $j^{-1}e^{-1}\zeta_*[N_0]=\partial x$ for some $x\in H_5(C,\partial)$.
By homology Alexander duality $x=kA[N]$ for some integer $k$.
 We have $k=1$ because
$$k[N]=\nu_*\partial(kA[N])=\nu_*j^{-1}e^{-1}\zeta_*[N_0]=
%j_N^{-1}e_N^{-1}
\nu|_{N_0,*}\zeta_*[N_0]=
%j_N^{-1}e_N^{-1}
[N_0]=[N].$$
\quad
Now we prove (b).
Observe that the normal bundle $\nu_\zeta$ of embedding
$\zeta:N_0\to\partial C$ is isomorphic to $\zeta^\perp$.
Hence their homology Euler classes coincide.
Then by (a) we have $\zeta^!ej\partial A[N]=\zeta^!\zeta_*[N_0]=PDe(\nu_\zeta)$.
Also,
$$\varkappa (f)\overset{(1)}\to=\nu_*(\partial A[N])^2\overset{(2)}\to=
\nu|_{N_0,*}ej(\partial A[N])^2
\overset{(3)}\to=\nu|_{N_0,*}(ej\partial A[N])^2\overset{(4)}\to=
\nu|_{N_0,*}(\zeta_*[N_0])^2\overset{(5)}\to=PDe(\nu_\zeta).$$
Here

$\bullet$ squares denote the intersection squares;

$\bullet$ the first equality is the definition of $\varkappa $;

$\bullet$ the second equality holds because $\nu_*=\nu|_{N_0,*}ej$;

$\bullet$ the third equality holds by the naturality properties of $\cap$;
%\footnote{Here is a proof of $\partial(A[N])^2=(\partial A[N])^2$
%(analogously one checks that $ej(\partial A)^2=e(j\partial A)^2=(ej\partial A)^2$).
%Take a triangulation of $C$.
% Represent $A$ by a 5-chain $a$ in this triangulation and a 5-chain $a'$ in the dual triangulation.
%Then $(\partial A[N])^2$ is represented by $\partial a\cap\partial a'=\partial (a\cap a')$ which represents $\partial( A[N]^2)$.}

$\bullet$ the fourth equality follows by (a);

$\bullet$ the fifth equality is the definition of $PDe(\nu_\zeta)$.
\qed

%$$\minCDarrowwidth{0pt}\CD
%H_4(N_0,\partial) @>> \xi_* > H_4(\nu^{-1}N_0,\partial) @<< e <
%H_4(\partial C,\nu^{-1}B^4) @<< j <  H_4(\partial C) @>> >\\
%@VV \cong V    @VV  V    @VV V @VV V    @VV  V\\
%H_4(N)  @>> ej\overline\xi_* > H_4(C,\nu^{-1}B^4)  @<< e <
%H_4(C\cup\widehat\nu^{-1}B^4,\widehat\nu^{-1}B^4) @<< j <
%H_4(C\cup\widehat\nu^{-1}B^4)  @<< i_* < H_4(C)\endCD.$$
%Here $e$ and $j$ are isomorphisms; vertical arrows are induced by the
%inclusions.

\bigskip
{\bf Compressible embeddings.}

Assume that $H_1(N)=0$.
A (smooth) embedding $f:N\to\R^7$ is called {\it PL compressible} if for some
embedding $g:S^4\to\R^7$ the embedding $f\#g$ is isotopic to an embedding
$f':N\to\R^7$ such that $f'(N)\subset\R^6$ (this is equivalent to saying that
$f$ is PL isotopic to such an embedding $f'$).
%We need $H_1(N)=0$ because there could be an obstruction in $C^3_3$
Cf. the Effectiveness Theorem \Eff.
The study of compressible embeddings is a classical problem in topology of
manifolds, see references in [Sk08'].
%[Ha66, Hi66, Gi67, Ti69, Vr89, RS01, CR05, Ta06, .

\smallskip
{\bf Remark \Cor.} [Vr89]
{\it Let $N$ be a closed connected 4-manifold such that $H_1(N)=0$.

(a) An embedding $f:N\to\R^7$ is PL compressible if and only if $\varkappa (f)=0$
(which holds if and only if $\pi_3(C_f)\cong\Z$).

(b) Two PL compressible embeddings $f:N\to\R^7$ differ only by connected sum
with an embedding $S^4\to\R^7$.

(c) The map $ E^6_{PL}(N)\to E^7_{PL}(N)$ induced by the inclusion
$\R^6\to\R^7$ is trivial.}

\smallskip
{\it Proof.}
Part (b) follows from part (a) and the Bo\'echat-Haefliger Theorem \Bh(b).
Part (c) follows from the Bo\'echat-Haefliger Theorem \Bh(b) and the PL version
of (a), which is proved analogously to (a).
Let us prove (a).

By the Complement Lemma \Com(b) $\varkappa (f)=0$ is equivalent to $\pi_3(C_f)=\Z$.

Clearly, for an embedding $f:N\to S^7$ isotopic to $f'\#g$ for some embeddings
$g:S^4\to\R^7$ and $f':N\to\R^7$ such that $f'(N)\subset\R^6$, we have
$C_f\simeq\Sigma(S^6-f'N)$.
Hence $\varkappa (f)=\varkappa (f')=0$ by the Complement Lemma \Com(a).

If $\varkappa (f)=0$, then by the Bo\'echat-Haefliger Theorem \Bh(a) we have
$w_2(N)=0$ and $\sigma(N)=0$.
Hence there is an embedding $f':N\to S^6$ [CS79, Ru82, GS99, Remark 9.1.22].
We have $\varkappa ((i:\R^6\to\R^7)\circ f')=0$.
Hence by the Bo\'echat-Haefliger Theorem \Bh(b) $f$ is PL compressible.
\qed

%\newpage
\head 3. Attaching invariant and proof of the Effectiveness Theorem \Eff \endhead

Recall that {\it $N$ is a closed connected orientable 4-manifold and
$f:N\to\R^7$ is an embedding}.
Fix an orientation of $N$ and of $\R^7$.
Take a small oriented disk $D^3_f\subset\R^7$ whose intersection
with $fN$ consists of exactly one point of sign $+1$ and such that
$\partial D^3_f\subset\partial C_f$.
{\it The meridian} $S^2_f$ of $f$ is $\partial D^3_f$.
Identify $S^2_f$ with $S^2$.

Let $G_q$ be the space of maps $S^{q-1}\to S^{q-1}$ of degree +1.
The space $G_q$ is identified with a subspace of $G_{q+1}$ via suspension.
Let $G=\lim_q G_q$ and $SO=\lim_q SO_q$.
% and $G/SO=\cup_q G_q/SO_q$.
The base points are the identity map or its class.

By $*$ we denote the base point of any space.
By $[X,Y]$ we denote the set of {\it based} homotopy classes of maps $X\to Y$
(the choice of base points is clear from the context).

%$\pi_i(G_q/SO_q\simeq G_q\cup\Con SO_q)\ne\pi_i(G_q,SO_q)$ !!!

\bigskip
{\bf The Haefliger-Levine attaching invariant of knots.}

\smallskip
{\it Construction of the attaching invariant $a: E^7(S^4)\to\pi_4(G_3,SO_3)$.}
Take a smooth embedding $f:S^4\to S^7$.
The space $C_f$ is simply-connected and by Alexander duality the inclusion
$S^2_f\to C_f$ induces an isomorphism in homology.
Hence this inclusion is a homotopy equivalence.
Take a homotopy equivalence $h_f$ homotopy inverse to the inclusion
$S^2_f\to C_f$.
We may assume that $h_f$ is a retraction onto $S^2_f$.
Since orientations of $S^7$ and of $S^4$ are fixed, the homotopy
class of $h_f$ depends only on $f$.
Since $\nu_f$ is trivial [Ke59, Ma59], there is a framing
$\xi:S^4\times S^2\to\partial C_f$ of $\nu_f$.
We may assume that $\xi(*\times S^2)=S^2_f$.

%Since $e(\nu_f)=0\in H_1(S^4)$, there is a section of $\nu_f$ wrong!!!
%Since $H^{i+1}(S^4;\pi_i(S^1))=0$, it follows that a section of $\nu_f$ can be
%extended to a framing (in particular, that

{\it The attaching invariant} $a(f,\xi)$ is the homotopy class of the
composition
$$S^4\times S^2 \overset\xi\to\cong\partial C_f\subset C_f\overset{h_f}
\to\simeq S^2_f.$$
Clearly, $a(f,\xi)$ is independent on isotopy of $f$.
Since $a(f,\xi)|_{*\times S^2}=\id S^2$, the map $a(f,\xi)|_{x\times S^2}$ is a
homotopy equivalence of degree +1 for each $x$.
Hence $a(f,\xi)\in\pi_4(G_3)$.

The choice of a framing $\xi$ is in $\pi_4(SO_3)$.
Since the composition $\pi_4(SO_3)\to\pi_4(G_3)\to\pi_4(G_3,SO_3)$ is
trivial, it follows that the image $a(f)\in\pi_4(G_3,SO_3)$ of $a(f,\xi)$
does not depend on $\xi$.
This class $a(f)$ is called {\it the attaching invariant} of $f$.

\smallskip
Clearly, $a: E^7(S^4)\to\pi_4(G_3,SO_3)$ is a homomorphism.
The injectivity of $a$ is proved in \S4.
The surjectivity of $a$ can be proved analogously, cf. [Fu94].
%Cf. the Symmetry Remark in \S3.
We prove the following well-known lemma because we could not find the proof in
the literature.

\smallskip
{\bf Lemma \Gso.} \quad $\pi_5(G,SO)=0$, \quad $\pi_4(G,SO)\cong\Z$ \ and
\ $\pi_4(G_3,SO_3)\cong\pi_6(S^2)\cong\Z_{12}$, \quad
cf. [Ha66, the text before Corollary 6.6].

\smallskip
{\it Proof.}
Let $F_q$ be the space of maps $S^q\to S^q$ of degree +1 leaving the north pole fixed.
From the homotopy exact sequence of fibration $F_q\to G_{q+1}\to S^q$ we get
%obtain
$\pi_n(G_{q+1})\cong\pi_n(F_q)$
%\cong\pi_{n+q}(S^q)$
for $q>n+1$.
Since $F_q\simeq \Omega_qS^q$, we have $\pi_n(F_q)\cong\pi_{n+q}(S^q)$ for $n>0$.
So $\pi_n(G)=\pi_n^S$.

%In order to prove that $\pi_5(G,SO)=0$ and $\pi_4(G,SO)\cong\Z$
Consider the exact sequence of pair
$$\minCDarrowwidth{0pt}\CD
\pi_5(G) @>> > \pi_5(G,SO) @>> >\pi_4(SO) @>> > \pi_4(G) @>> > \pi_4(G,SO)
@>> > \pi_3(SO) @>> > \pi_3(G)\\
@VV \cong V    @.          @VV \cong V    @VV \cong V       @.
@VV \cong V    @VV \cong V\\
\pi_5^S=0       @.        @.        0    @.  \pi_4^S=0  @. @. \Z
@. \pi_3^S\cong\Z_{24}
\endCD$$
We obtain that $\pi_5(G,SO)=0$ and $\pi_4(G,SO)$ is isomorphic to a
subgroup of $\pi_3(SO)\cong\Z$ having a finite index, i.e. to $\Z$.

%In order to prove that $\pi_4(G_3,SO_3)\cong\pi_6(S^2)$
Consider the fibration
$F_2\to G_3\to S^2$ and its subfibration $SO_2\to SO_3\to S^2$.
We obtain the following diagram
$$\minCDarrowwidth{0pt}\CD
@. \pi_3(SO_2) @. \pi_3(SO_3) @. \\
@. @A 0 AA         @A AA          @. \\
0 @>> > \pi_4(F_2,SO_2) @>> \cong > \pi_4(G_3,SO_3) @>> > 0\\
@A AA @A \cong AA @A AA @A AA \\
\pi_5(S^2) @>> \partial > \pi_4(F_2) @>> i > \pi_4(G_3) @>> p > \pi_4(S^2) \\
@A = AA @A 0 AA @A AA @A = AA \\
\pi_5(S^2) @>> > \pi_4(SO_2) @>> > \pi_4(SO_3) @>> > \pi_4(S^2) \\
\endCD$$
Since $\pi_i(SO_2)=0$ for $i\ge2$,
$\pi_4(G_3,SO_3)\cong\pi_4(F_2,SO_2)\cong\pi_4(F_2)\cong\pi_6(S^2)$.
%The last isomorphism follows analogously to the first paragraph of the
%previous proof.
\qed

\smallskip
{\bf Symmetry Remark.}
{\it For each embedding $g:S^4\to\R^7$ the composition of $g$ with the
reflection-symmetry $\R^7\to\R^7$ is isotopic to $-g$.
\footnote{By definition, $-g$ is the composition of $g$ with
reflection-symmetries of $S^4$ and of $\R^7$ [Ha66].
So change of the orientation on $S^4$ alone does not change
the isotopy class of $g$.}
}

\smallskip
{\it Proof.}
%We use the definition of the attaching invariant
%$a:E^7(S^4)\to\pi_4(\overline{G_3})$ from this section???. the following subsection?
Recall that the attaching invariant $a$ is an isomorphism.
Change of orientation of $\R^7$ induces change of orientation of $S^2_f$.
Since the
%unlinked
framing respects orientations, under change of orientation
of $\R^7$ the
%unlinked
framing changes orientation at each point.
Therefore change of orientation of $S^7$ carries the attaching invariant $a=a(\psi):S^4\times S^2\to S^2$
to $\sigma_2\circ a\circ(\id S^4\times\sigma_2)$, where $\sigma_m:S^m\to S^m$ is the reflection-symmetry.
Identify $\pi_4(G_3,SO_3)$ with $\pi_6(S^2)$ by the isomorphism from Lemma \Gso.
Under this identification $\id S^4\times\sigma_2$ goes to $\sigma_6$.
Then the required relation follows because
$$\sigma_2\circ a\circ\sigma_6=\sigma_2\circ(-a)=a-[\iota_2,\iota_2]\circ h_0(a)=
a-2\eta\circ h_0(a)=a-2a=-a\in\pi_6(S^2).$$
Here $h_0:\pi_6(S^2)\to\pi_6(S^3)$ is the generalized Hopf invariant, which is an isomorphism
inverse to the composition with the Hopf map $\eta\in\pi_3(S^2)$
[Po85, Complement to Lecture 6, (10)].
\qed

%{\it An alternative proof that $\pi_n(G)\cong\pi_n^S$.}
%Consider the Barratt-Puppe exact sequence of the pair
%$(S^n\times S^{q-1},S^n\times*)$ and maps to $S^{q-1}$ (with base points).
%Recall that $S^n\times S^{q-1}/S^n\times*\simeq S^{n+q-1}\vee S^{q-1}$.
%Since $\pi_n(S^{q-1})=\pi_{n+1}(S^{q-1})=0$ for large $q$, it follows that
%the action $[S^{n+q-1}\vee S^{q-1};S^{q-1}]\to[S^n\times S^{q-1};S^{q-1}]$
%from the exact sequence is a 1--1 correspondence.
%Hence homotopy classes of maps
%$S^n\times S^{q-1}\to S^{q-1}$ identical on
%$*\times S^{q-1}$ are in 1--1 correspondence with $\pi_{n+q-1}(S^{q-1})$.
%%(is a homomorphism between cohomotopy groups in the sence of Borsuk for
%%$q\ge n+2$),
%This 1--1 correspondence gives an isomorphism
%$\pi_n(G)\cong\pi_{n+q-1}(S^{q-1})\cong\pi_n^S$.

\bigskip
{\bf Proof of the Effectiveness Theorem \Eff.}

For $X\subset N$ an orientation-preserving normal framing
$\xi:X\times S^2\to\partial C_f$ is called {\it unlinked} if the composition of
the section $\xi_1:X\to\partial C_f$ (formed by first vectors of the framing)
with the inclusion $\partial C_f\subset C_f$ is null-homotopic.

%\footnote{This is in some sense analogous to {\it unlinked} normal vector field and framing on a knot in $S^3$.}

The following two lemmas are proved using standard arguments.

\smallskip
{\bf Extension Lemma \Exn.}
{\it (a) Let $g:S^4\to\R^7$ and $f:N\to\R^7$ be embeddings such that
$f(N)\subset\R^6$.
Then there is an unlinked framing $\xi:N_0\times S^2\to\partial C_{f\#g}$.

(b) If $\sigma(N)=0$, then any framing of $\nu_f|_{N_0}$ extends to that of
$\nu_f$.}

\smallskip
{\it Proof.}
(a) Any embedding $N\to\R^6$ has trivial normal bundle.
Thus $f$ has a framing $\xi:N_0\times S^2\to\partial C_f$ such that the
section $\xi_1:N_0\to\partial C_f$ formed by third vectors is orthogonal
to $\R^6$.
Then the composition of $\xi_1$ with the inclusion $\partial C_f\subset C_f$ is
null-homotopic.
We may assume that $(f\#g)|_{N_0}=f|_{N_0}$ and $(f\#g)(N-N_0)$ misses the
trace of the null-homotopy.
Hence $\xi_1(N_0)\subset \partial C_f\cap \partial C_{f\#g}$ and
$\xi:N_0\times S^2\to\partial C_{f\#g}$ is an unlinked framing.
\qed

%\smallskip
%{\it Proof of (b).}
(b) Given a framing $\xi$ of $\nu_f|_{N_0}$, there is a complete obstruction
$O=O(\xi)\in H^4(N,\pi_3(SO_3))\cong\Z$ to extension of $\xi$ to $N$.
Since the inclusion $SO_3\subset SO$ induces a multiplication by 2 on $\pi_3$,
it follows that $O$ equals to twice the obstruction to extension of $\xi$ to a
{\it stable} framing of $\nu_f$.
Hence $\pm O=p_1(N)=3\sigma(N)=0$ [Ma80, argument before Lemma 1.15].
\qed

\smallskip
{\bf The Barratt-Puppe Sequence Lemma \Bp.}
{\it Let $P$ and $B$ be CW complexes with base points $*$,
$$Y:=P\times B,\quad \Con B=B\times[0,1]/B\times1\quad\text{and}
\quad Y/B:=Y\bigcup\limits_{*\times B=B\times0}\Con B.$$
Then the following sequence of sets is exact:
$$[\Sigma B\times P;P]_P\overset\rho\to\to[\Sigma B;P]\overset\psi\to\to [Y/B;P]_P\overset R\to\to [Y;P]_P, \quad\text{where}$$
\quad
$\bullet$ subscript $P$ means that we consider retractions to $P=[\{-1\}\times B]\times P\subset \Sigma B\times P$, to $P=*\times P\subset Y/B$ and to $P=*\times P\subset Y$, respectively, up to homotopy fixed on $P$;

$\bullet$ the base point of $[\Sigma B;P]$ is the constant map to $*$.

$\bullet$ the base points of $[Y;P]_P$, of $[Y/B;P]_P$ and of $[\Sigma B\times P;P]_P$ are
the projection $p:Y\to P$, the extension $p':Y/B\to P$ of $p$ that maps $\Con B$ to $*$, and the projection to $P$, respectively;

$\bullet$ $\rho$ and $R$ are the restriction maps;

$\bullet$ $\psi(f)$ is the composition $Y/B\to (Y/B)\vee \Sigma B\overset{p'\vee f}\to\to P$ of the contraction of $B\times \frac12\subset\Con B$ and of $p'\vee f$. }

\smallskip
{\it Proof.} Denote
$$\Con\phantom{}_-B=B\times[-1,0]/B\times\{-1\}\quad\text{and}\quad
Y':=\Con\phantom{}_-B\times P\cup\Con B.$$
We can replace $[\Sigma B;P]$ by $[\Sigma B\vee P;P]_P$ and further by $[Y';P]_P$, where $P=[\{-1\}\times B]\times P.$
Then $\rho$ and $\psi$ will be replaced by restriction maps, which will be denoted by the same letters $\rho$ and $\psi$.

The restriction to $Y$ of a retraction $f:Y/B\to P$ is homotopic to $p$ if and only if $f$ extends to
a retraction $Y'\to P\times\{-1\}$.
Hence the obtained sequence is exact at $[Y/B;P]_P$.

The restriction to $Y/B$ of a retraction $f:Y'\to P\times\{-1\}$ is homotopic to $p'$ if and only if $f$ extends to a retraction $\widetilde Y:=\Sigma B\times P\bigcup\limits_{\Con B\times *=B*2}B*[2,3]$.
Since $\widetilde Y$ deformation retracts to $\Sigma B\times P$ relative to $P$, the obtained sequence is exact at $[Y';P]_P$.
QED
\footnote{In the previous version of this paper I used the version of the Barratt-Puppe Sequence Lemma \bp with
$[\Sigma B\times P;P]_P$ replaced by $[\Sigma Y;P]$ (the base point is the constant map).
That version is false: the obtained sequence of based sets is not necessarily exact at $[\Sigma B;P]$.
(Note that the standard argument shows that the new sequence is exact at $[Y/B;P]_P$.)
Let me present a counterexample by D. Crowley.
Take $P=S^2$ and $B=S^1$.
Then the new sequence is
$$[\Sigma(S^2\times S^1);S^2]\overset\rho\to\to[\Sigma S^1;S^2]\overset\psi\to\to [S^3\vee S^2;S^2]_{S^2}
\overset R\to\to [S^2\times S^1;S^2]_{S^2}.$$
(Note that $R$ is change on the top cell.)
Since $\Sigma(S^2\times S^1)$ retracts to $\Sigma S^1$, the map $\rho$ is surjective.
The final term is isomorphic to $\pi_1(G_3)\cong\Z_2$.
The last term but one is isomorphic to $\pi_3(S^2)\cong\Z$.
Thus $R$ is not injective.
Since the sequence is exact at $[S^3\vee S^2;S^2]_{S^2}$, the map $\psi$ is not constant.
So the sequence is not exact at $[\Sigma S^1;S^2]$, even the map $\psi\circ\rho$ is not constant.}

\smallskip
{\bf Retraction Lemma \Ret.} {\it Let $N$ be a closed simply-connected 4-manifold,
$f:N\to\R^7$ an embedding and $\xi:N_0\times S^2\to\partial C_f$ an unlinked framing.
 Then there is a unique (up to homotopy fixed on $S^2_f$) retraction
$r(\xi):C_f\to S^2_f$ whose restriction to $\xi(N_0\times S^2)$ is
the projection to $\xi(*\times S^2)=S^2_f$.
\footnote{The
unlinkedness of $\xi$ is essential in the Retraction Lemma \Ret.
Indeed, for an embedding $f:S^2\times S^2\to\R^7$ such that
$\varkappa (f)\ne0$ there is a framing of $\nu_f$ over $N_0$, however $C_f$
does not retract to $S^2$ because $C_f\simeq C_{\varkappa (f)}$  by the
Complement Lemma \Com(a).
\newline
The homotopy between retractions from the Retraction Lemma \ret is not
assumed to be fixed on $\xi(N_0\times S^2)$.
(The set of retractions $C_f\to S^2$ extending the projection
$N_0\times S^2\to S^2$ up to homotopy fixed on $N_0\times S^2$
is in 1--1 correspondence with $H^3(C_f,N_0\times S^2)\cong H_2(N)$.)
\newline
In the previous version of this paper the Retraction Lemma has the orientability assumption instead of
the simply connectedness.
The proof was incorrect, see the previous footnote.
The statement was also false [CS].}
}

\smallskip
{\it Proof.}
Denote $A:=\xi(N_0\times S^2)$ and identify $N_0$ with $\xi(N_0\times*)$.
Since the framing $\xi$ is unlinked, the inclusion $A\to C_f$ extends to a map $A\cup\Con N_0\to C_f.$
By the Alexander duality and the Mayer-Vietoris sequence this map induces a homology isomorphism.
Hence by the relative Hurewicz Theorem this map is a homotopy equivalence.

Since projection $p:A\to S^2$ is null-homotopic on $N_0$, map $p$ extends to $A\cup\Con N_0$.
This implies the existence of $r(\xi)$.
In order to prove the uniqueness consider the exact sequence of sets given by the Barratt-Puppe
Sequence Lemma \Bp:
$$[\Sigma N_0\times S^2;S^2]_{S^2}\overset\rho\to\to[\Sigma N_0;S^2]\overset\psi\to\to [A\cup\Con N_0;S^2]_{S^2}\overset R\to\to [A;S^2]_{S^2}.$$
Since $N$ is simply-connected, $N_0$ is homotopy equivalent to a wedge of 2-spheres.
Since $\Sigma[\eta,\iota_2]=0\in\pi_5(S^3)$ and $\Sigma:\pi_4(S^2)\to\pi_5(S^3)$ is an isomorphism [To62], we have
$[\eta,\iota_2]=0\in\pi_4(S^2)$.
So $[\alpha,\iota_2]=0\in\pi_4(S^2)$ for each $\alpha\in\pi_3(S^2)=\left<\eta\right>$.
Hence each map $\alpha\vee\id_{S^2}:S^3\vee S^2$ extends up to homotopy to $S^3\times S^2$.
Thus each map $\alpha\vee\id_{S^2}:\Sigma N_0\vee S^2$ extends up to homotopy to $\Sigma N_0\times S^2$, i.e., $\rho$ is surjective.
Therefore by exactness $R^{-1}(p)=p'$, which proves the uniqueness.
\qed

\smallskip
%\newpage
{\it Definition of the attaching invariant for an embedding $f:N\to\R^7$
which has an unlinked framing $N_0\times S^2\to\partial C_f$.}
Extend the unlinked framing of $f$ to a framing $\xi$ of $f$ by the Extension
Lemma \Exn(b).
We may assume that the $\nu_f$-preimage of the base point $*\in N_0$ is $S^2_f$.
Take the retraction $r=r(\xi|_{N_0})$ given by the Retraction Lemma \Ret.
{\it The attaching invariant} $a(f,\xi)$ is the homotopy class of the
composition
$$N\times S^2\overset\xi\to\cong\partial C_f\subset C_f\overset{r(\xi)}
\to\to S^2_f.$$
Since $r(\xi)$ is a retraction, $a(f,\xi)|_{*\times S^2}=\id S^2$.
Hence $a(f,\xi)|_{x\times S^2}$ is a homotopy equivalence
of degree +1 for each $x$.
Thus any map representing $a(f,\xi)$ can be identified with base point
preserving map $N\to G_3$.
Since $a(f,\xi)|_{*\times S^2}=\id S^2$ throughout a homotopy of $r(\xi)$
fixed on $S^2_f$, we may assume that $a(f,\xi)\in[N,G_3]$.

The choice of a framing $\xi$ is in  $[N,SO_3]$.
Under a change $\varphi:N\to SO_3$ a map $a:N\to G_3$ changes to the map
$a^\varphi:N\to G_3$ defined by $a^\varphi(x)=\varphi(x)a(x)$.
Let $\overline G_3$ be the homotopy fiber of $BSO_3\to BG_3$.
Thus $\pi_i(\overline G_3)\cong\pi_i(G_3,SO_3)$ and more generally there is an
exact sequence $[N,\Omega BSO_3]\to[N,\Omega BG_3]\to [N,\overline G_3]$.
Since $\Omega BG_3\simeq G_3$ and $\Omega BSO_3\simeq SO_3$, we may assume that
$a(f)\in [N,\Omega BG_3]$ and the choice of a framing $\xi$ is in
$[N,\Omega BSO_3]$.
Hence the image $a(f)\in[N,\overline G_3]$ of $a(f,\xi)$ does not depend
on $\xi$.
This class $a(f)=a_N(f)$ is called {\it the attaching invariant} of $f$.

\smallskip
{\bf Additivity Lemma \Add.} {\it For an embedding $f:N\to\R^7$ which has an
unlinked framing $N_0\times S^2\to\partial C_f$ and an embedding $g:S^4\to S^7$
we have $a_N(f\#g)=a_{S^4}(g)\#a_N(f)$, where
$\#:\pi_4(\overline G_3)\times[N,\overline G_3]\to[N,\overline G_3]$ is the
action given by the map $N\to N/\partial B^4\simeq N\vee S^4$.}

\smallskip
{\it Proof.} It follows because by definition of $a(f,\xi)$ we have
$a([f,\xi]\#[g,\zeta])=a(f,\xi)\#a(g,\zeta)$, where
$\#:\pi_4(G_3)\times[N,G_3]\to[N,G_3]$
is the action given by the map $N\to N/\partial B^4\simeq N\vee S^4$.
\qed

\smallskip
{\it Proof of the Effectiveness Theorem \Eff.} Let $f:N\to\R^7$ be
an embedding isotopic to $f_0\#g$ for some embedding $g:S^4\to\R^7$.
Take a framing given by the Extension Lemma \Exn(a). Then attaching
invariants $a(f,\xi)$ and $a(f)$ are defined.
We have $N/N_0\cong S^4$.
Consider the Barratt-Puppe exact sequence of sets of based homotopy classes:
$$[\Sigma N;\overline G_3]\overset{\Sigma R}\to\to[\Sigma N_0;\overline G_3]\to
\pi_4(\overline G_3)\overset v\to\to[N;\overline G_3]\overset R\to\to[N_0;\overline G_3].$$
By the Retraction Lemma \ret $a(f,\xi)|_{N_0\times S^2}$ is homotopic to the
projection onto $S^2$.
Thus the image of $a(f,\xi)$ under the restriction-induced map
$[N,G_3]\to[N_0,G_3]$ is the trivial homotopy class.
So $a_N(f)\in R^{-1}(*)=\im v$.

Since $N$ embeds into $\R^6$, we have $w_2(N)=0$.
Since $N$ is also simply-connected, by [Mi58] $\Sigma N$ retracts to $\Sigma N_0$.
Hence $\Sigma R$ is surjective.
So by exactness $v^{-1}(*)=0$.
Since $v$ extends to an action of the domain on the range, $v$ is injective.
The map $a_{S^4}: E^7(S^4)\to\pi_4(\overline G_3)$ is a monomorphism (this is
proved in [Ha66] or at the beginning of \S4).
So the theorem follows by the Additivity Lemma \Add.
\qed
\footnote{This proof shows that we can weaken the condition
{\it `$N$ embeds into $\R^6$ and $fN\subset\R^6$'} to
{\it `$fN_0\subset\R^6$ with trivial normal bundle and $\sigma(N)=0$'}.
Since the triviality of the normal bundle implies that $w_2(N)=0$, by the
remarks in \S1 the new assumption still implies that $N$ embeds into $S^6$.}

\bigskip
%\newpage
{\bf Higher-dimensional homology spheres.}

\smallskip
{\bf Theorem \Ntr.}
{\it Let $N$ be a homology $n$-sphere, $n\ge3$ and suppose that
$\Sigma^\infty:\pi_{n+2}(S^2)\to\pi_n^S$ is not injective.
Then for any embedding $f:N\to\R^{n+3}$ there is an embedding
$g:S^n\to\R^{n+3}$ such that $f\#g$ is not isotopic to $f$.}

\smallskip
{\it Proof.}
Analogously to the case $N=S^4$ we have that $\nu_f$ is trivial and $C_f\simeq S^2$
for each embedding $f:N\to\R^{n+3}$.
So analogously to the construction of the attaching invariant for $S^4$ we construct
$h_f:C_f\to S^2_f$ and take a framing $\xi$ of $\nu_f$.
 Since $H^i(N_0)=0$ for $i>0$,
%need twisted?
every two maps $N_0\times S^2\to S^2$ are homotopic.
In particular, a homotopy equivalence $h_f:C_f\to S^2$ is homotopic on $\xi(N_0\times S^2)$
to the projection onto $S^2$.
Now we construct the attaching invariant $a(f,\xi)\in[N,G_3]$ analogously to $S^4$.
 \footnote{Note that by Homotopy Lemma of \S5 and the following lemma we can obtain an attaching invariant in $\pi_{n+2}(S^2)$ directly.
\newline
{\bf Lemma.} {\it For $n\ge3$ and each embedding $f:N\to S^{n+3}$ of an $n$-dimensional homology sphere $N$
there is a unique unlinked framing $\xi_f:N\times S^2\to\partial C_f$.}
\newline
{\it Proof.}
By [Ke59, Ma59] $\nu_f$ is trivial.
The difference elements for sections $N\to\partial C_f$ are in
$H^i(N,\pi_i(S^2))$.
Hence the sections are in 1-1 correspondence with elements of $\pi_n(S^2)$.
Recall that $C_f\simeq S^2$.
So there is a unique section $\xi_{1,f}:N\to\partial C_f$ whose composition
with the inclusion $\partial C_f\subset C_f$ is null-homotopic.
Since $H^{i+1}(N,\pi_i(S^1))=0$ and $H^i(N,\pi_i(S^1))=0$, this section can be
uniquely extended to an unlinked framing $\xi_f:N\times S^2\to\partial C_f$.
\qed
}

Consider the Barratt-Puppe exact sequence of sets of based homotopy classes:
$$[\Sigma N_0;G_3]\to\pi_4(G_3)\overset v\to\to[N;G_3]\overset R\to\to[N_0;G_3].$$
Since $a(f,\xi)|_{N_0\times S^2}$ is homotopic to the projection onto $S^2$,
the image of $Ra(f,\xi)=*$.
So $a(f,\xi)\in R^{-1}(*)=\im v$.
Since $N$ is a homology sphere, $\Sigma N_0$ is contractible.
So by exactness $v^{-1}(*)=0$.
Since $v$ extends to an action of the domain on the range, $v$ is injective.
Let $a_N(f)\in\pi_n(G_3,SO_3)$ be the image of $v^{-1}a(f,\xi)$.

Clearly, $a_N(f\#g)=a_N(f)+a_{S^4}(g)$.
%Now (instead of using the injectivity of $a_{S^4}$) we use
Consider the exact sequence [Ha66, 4.11, cf. Ha86]
$$\pi_{n+1}(G,SO)\to  E^{n+3}(S^n)\overset{a_{S^n}}\to\to\pi_n(G_3,SO_3)\overset{st}\to\to\pi_n(G,SO),$$
where $st$ is the stabilization map.
Recall that $st$ equals to a composition
$$\pi_n(G_3,SO_3)\overset\cong\to\to\pi_{n+2}(S^2)\overset{\Sigma^\infty}\to\to \pi^S_n\to\pi_n(G,SO)$$
(Indeed, in the following commutative diagram
$$\minCDarrowwidth{2pt}\CD
\pi_{n+2}(S^2) @>>\cong> \pi_n(F_2) @>>i> \pi_n(G_3) @>>j> \pi_n(G_3,SO_3)\\
@VV\Sigma^\infty V      @VV st_F V        @VV st_G V          @VV st V        \\
\pi_n^S       @>>\cong> \pi_n(F) @>>i> \pi_n(G) @>>j> \pi_n(G,SO) \endCD$$
the upper map $ji$ is an isomorphism analogously to the end of the proof of Lemma \Gso.)

Since $\Sigma^\infty$ is not monomorphic, $st$ is not monomorphic.
Thus $\im a_{S^n}\ne0$.
This implies the Theorem.
\qed

\smallskip
The assumption of Theorem \ntr holds for each $n\le19$ except
$n\in\{6,7,9,15\}$ [To62, tables].
For $n=3$ Theorem \ntr is covered by [Ha72, Ta06].
Analogously to Theorem \ntr it follows that
{\it for any homology $n$-sphere $N$ the action of $E^{n+3}(S^n)$ on
$E^{n+3}(N)$ is}

{\it (a) non-trivial, if the stabilization map
$\pi_n(\overline G_3)\to\pi_{n+1}(G,SO)$ is not injective.

(b) effective, if $\pi_{n+1}(G,SO)=0$.}

\comment

\smallskip
{\it Proof (suggested by a referee).}
Denote $A:=\xi(N_0\times S^2)$.
Since the framing $\xi$ is unlinked, the inclusion $A\to C_f$ extends to
a map $A\cup\Con(N_0\times*)\to C_f.$
By the Alexander duality and the Mayer-Vietoris sequence
this map induces a homology isomorphism.
Hence by the relative Hurewicz Theorem this map is a homotopy equivalence.

Since projection $p:A\to S^2$ is null-homotopic on $N_0\times*$, $p$
extends to $A\cup\Con(N_0\times*)$. This implies the existence in
the Retraction Lemma \Ret.

In order to prove the uniqueness
consider the Barratt-Puppe exact sequence of sets???:
$$[\Sigma A;S^2]\overset\rho\to\to[\Sigma(N_0\times*);S^2]\to
[A\cup\Con(N_0\times*);S^2]_{S^2}\overset R\to\to [A;S^2]_{S^2}.$$
Here $\rho$ and $R$ are the restrictions; the subscript $S^2$ means
that we consider retractions to $S^2=*\times S^2\subset A$ up to
homotopy fixed on $S^2$; the base points of $[A;S^2]_{S^2}$ and of
$[A\cup\Con(N_0\times*);S^2]_{S^2}$ are $p$ and any retraction $r$
whose existence is proved in the previous paragraph. Since $\Sigma A$ retracts to $\Sigma(N_0\times*)$,
\footnote{We have
$\Sigma(X\times Y)$ retracts to $\Sigma X$ for CW complexes $X$ and
$Y$. This follows because $\Sigma(X\times Y)\simeq \Sigma
X\vee\Sigma Y\vee\Sigma(X\wedge Y)$. The following arguments for
this were communicated by J. Klein.
\newline
{\it First proof.} Let $t:X\times Y\to X\wedge Y$ be the quotient map.
Let $p$ be a composition $\Sigma(X\times Y)\to\Sigma X\to \Sigma X\vee\Sigma Y$ of the projection
(defined by $(x,y,s)\mapsto(x,s)$) and the inclusion.
Let $q:\Sigma(X\times Y)\to\Sigma Y\to\Sigma X\vee\Sigma Y$ be analogous composition.
Then $(p+q)\vee\Sigma t:\Sigma(X\times Y)\to\Sigma X\vee\Sigma Y\vee\Sigma(X\wedge Y)$ is a homology isomorphism.
By the relative Hurewicz theorem it is also a homotopy isomorphism.
Then by the Whitehead theorem it is a homotopy equivalence.
\newline
{\it Second proof.} Taking the Hopf construction of  the
identity map of $X\times Y$ we get a map $h:X * Y \to \Sigma(X\times Y)$.
The composition $X * Y\to \Sigma(X\times Y)\to \Sigma(X \wedge Y)$ of
$h$ and the quotient map is given by collapsing out
a suitable contractible subspace of the join.
So this composition is a homotopy equivalence.
The wedge $\Sigma X\vee\Sigma Y\vee (X*Y)\to\Sigma (X\times Y)$ of the inclusions and $h$ is a homology isomorphism.
Then, as above, this wedge is a homotopy equivalence. }
%The argument for this is inductive.
%The cells of $\Sigma(X\times S^k)$ are the cell of $\Sigma S^k$, the cells of $\Sigma X$ and
%product cells. The product cells have attaching maps which are suspensions over Whitehead
%products. It is wrong that the attaching maps of product cells are in general Whitehead products!!!
%Since the latter suspensions are null-homotopic, the inductive step follows.
%it follows that
$\rho$ is surjective. So by exactness $R^{-1}(p)=r$,
%Since $v$ extends to an action [why the domain is a group? but no need of the action]
%of the domain on the range, $v$ is injective.
%This means that a map $C_f\simeq A\cup\Con(N_0\times*)\to S^2$ is
%uniquely defined (up to homotopy fixed on $S^2$) by its restriction to $A$.
%Since the inclusion $A\to C_f$ extends to a homotopy equivalence $A\cup\Con(N_0\times*)\to C_f$,
which proves the Lemma.  \qed

%the above argument without action is not on arxiv

%We proved that $(C_f,A)\simeq(A\cup\Con(N_0\times*),A)$; it would be
%interesting to know if $C_f\simeq \Sigma^2N_0\vee S^2$.

%A referee suggested that an alternative proof of the uniqueness 'is given by
%considering the fibration up to homotopy of mapping spaces
%$$map(\Sigma N_0,S^2)\to map(\Sigma N_0\vee S^2,S^2)\to map(N_0\times S^2,S^2)$$
%noting that the first of these maps is null-homotopic (because
%$\Sigma(N_0\times S^2)$ retracts onto $\Sigma N_0$).'

\endcomment

%\newpage
\head 4. Proof of the Primitivity Theorem \Pri \endhead

%\bigskip
{\bf Preliminary results.}

A map is called {\it $m$-connected} if it induces an isomorphism
on $\pi_i$ for $i<m$ and an epimorphism on $\pi_m$.

\smallskip
{\bf Almost Diffeomorphism Theorem \Adk.}
{\it Let $C_0$ and $C_1$ be compact simply-connected 7-manifolds such that
$H_3(C_0)=H_3(C_1)=0$ and $\varphi:\partial C_0\to\partial C_1$ a
diffeomorphism.
For some homotopy 7-sphere $\Sigma$ there is a diffeomorphism
$\varphi:C_0\to C_1\#\Sigma$ extending $\varphi$ if and only if there exist

$\bullet$ a fibration $p:B\to BO$ such that $\pi_1(B)=0$,

$\bullet$ a compact 8-submanifold $W\subset S^{18}$ such that
$\partial W=M_\varphi:=C_0\cup_\varphi(-C_1)$, and

$\bullet$ a lifting $\overline\nu:W\to B$ of the
%Gauss
classifying map $\nu:W\to BO$ of the normal bundle such that
$\overline\nu|_{C_0},\overline\nu|_{C_1}$ are 4-connected. }

\smallskip
This is a slight improvement of [Kr99, Theorem 3 and Remark in p. 730]
(in which remark for $q$ even and $\pi_1(B)=0$ we can take $k\ge q-1$).
At the end of this section we present the author's write-up of M. Kreck's
complete proof of the Almost Diffeomorphism Theorem \Adk.

%because 4-connected maps $\overline\nu|_{C_0},\overline\nu|_{C_1}$ are
%3-smoothings not necessarily 4-smoothings [Kr99, Definition (i) in p. 711]
%and [Kr99, p. 708, first lines].
%$\bullet$ in this situation the Kreck surgery obstruction is elementary if and
%only if the ordinary surgery obstruction $\sigma(W)$ is trivial;
%$\bullet$ one can make $\sigma(W)$ trivial by adding boundary of a homotopy
%sphere.
%A direct proof (which essentially proves [Kr99, Remark in p. 730]) is given
%here

Denote by $BO\left<m\right>$ the (unique up to homotopy equivalence)
$(m-1)$-connected space for which there exists a fibration
$p:BO\left<m\right>\to BO$ inducing an isomorphism on $\pi_i$ for $i\ge m$.
(There is a misprint in [Kr99], Definition of $k$-connected cover
on p. 712: $X\left<k\right>$ should read as $X\left<k+1\right>$.)

For a fibration $p:B\to BO$ denote by $\Omega_q(B)$ the group of bordism
classes of liftings $\overline\mu:Q\to B$ of the classifying map of stable
normal bundle $\mu:Q\to BO$, where $Q$ is a (non-fixed) $q$-submanifold of
$\R^{3q}$.
Two such liftings $\overline\mu:Q\to B$ and $\overline\mu':Q'\to B$ are called
{\it bordant} if there is a $(q+1)$-submanifold $W\subset\R^{3q+3}$ and a
lifting $M:W\to B$ of the classifying map $W\to BO$ of stable normal bundle
such that $\partial W=Q\sqcup Q'$ and
$M|_{\partial W}=\overline\mu\sqcup\overline\mu'$.
This should be denoted by $\Omega_q(p)$ not $\Omega_q(B)$ but no confusion
would appear.
(This group is the same as $\Omega_q(B,p^*t)$ in the notation of [Ko88].)

\smallskip
{\bf Reduction Lemma \Red.}
{\it
%Let $q$ be an even integer, $4\le n\le 2q-4$ and $N$ a closed connected
%$n$-manifold.
Embeddings $f,f':N\to\R^7$ are isotopic if for some orientation-preserving
bundle isomorphism $\varphi:\partial C_f\to\partial C_{f'}$ and some embedding
$M_\varphi:=C_f\cup_\varphi(-C_{f'})\to S^{16}$ there exist

$\bullet$ a space $C$,

$\bullet$ a map $h:M_\varphi\to C$ whose restrictions to $C_f$ and to $C_{f'}$
are 4-connected, and

$\bullet$ a lifting $l:M_\varphi\to BO\left<5\right>$ of the classifying map
$M_\varphi\to BO$ of the normal bundle such that
$$[h\times l]=0\in\Omega_7(C\times BO\left<5\right>)/i_*\theta_7,$$
where $\theta_7$ is the group of orientation-preserving diffeomorphism classes
of homotopy 7-spheres.}

\smallskip
This situation is explained by the following diagram.
$$\minCDarrowwidth{5pt}\CD
C\times BO\left<5\right> @>> \pr_2 > BO\left<5\right> \\
@AA h\times l A                   @VV p V          \\
M_\varphi=C_f\cup_\varphi(-C_{f'})     @>>    > BO     \endCD$$

%\smallskip
{\it Proof of the Reduction Lemma \Red.}
Denote $B:=C\times BO\left<5\right>$.
Since $[h\times l]=0$, it follows that there is a homotopy 7-sphere
$\Sigma'$ and a map $s:\Sigma'\to B$ such that $(h\times l)\#s$ is
null-bordant.
Take a
%simply-connected
null-bordism $\overline\nu:W\to B$ of $(h\times l)\#s$.
Since $h|_{C_f}$ and $h|_{C_{f'}}$ are 4-connected and $BO\left<5\right>$ is
4-connected, it follows that $\overline\nu|_{C_f}=(h\times l)|_{C_f}$ and
$\overline\nu|_{C_{f'}}=(h\times l)|_{C_{f'}}$ are 4-connected.
Therefore by the Almost Diffeomorphism Theorem \adk $\varphi$ extends to an
orientation-preserving diffeomorphism $C_f\cong C_{f'}\#\Sigma'\#\Sigma$ for
some homotopy sphere $\Sigma$.
The isomorphism $\varphi$ also extends to an orientation-preserving
diffeomorphism $S^7-\Int C_f\to S^7-\Int C_{f'}$.
Therefore $\Sigma'\#\Sigma\cong S^7\#\Sigma'\#\Sigma\cong S^7$.
Hence $f$ is isotopic to $f'$ by Lemma \Is.
\qed
\footnote{This argument shows that in Lemma \is the condition for isotopy
could be weakened to `there is an orientation-preserving bundle isomorphism
$\varphi:\partial C_f\to\partial C_{f'}$ which extends to an
orientation-preserving diffeomorphism $C_f\to C_{f'}\#\Sigma$ for some homotopy
$n$-sphere $\Sigma$.'}

\smallskip
The remaining lemmas of this subsection are
essentially known and/or are proved using standard arguments.

\smallskip
{\bf Fiber Lemma \Fib.} {\it
Let $F$ be the fiber of $p:BO\left<5\right>\to BO$.
Then}
$$\pi_i(F)=0\quad\text{for}\quad i\not\in\{1,3\},\quad\pi_1(F)\cong\Z_2
\quad\text{and}\quad\pi_3(F)\cong\Z.$$

%\smallskip
{\it Proof.}
Recall that $\pi_i(O)\cong\pi_{i+1}(BO)\cong\Z_2,\Z_2,0,\Z,0,0,0,\Z$
according to $i=0,1,2,3$, $4,5,6,7$.
From the homotopy exact sequence of the fibration $F\to BO\left<5\right>\to BO$
we obtain the assertion for $i\ge4$. From the same exact sequence we obtain
$\pi_i(F)\cong\pi_{i+1}(BO)\cong\pi_i(SO)\cong\Z_2,0,\Z$ for $i=1,2,3$.

\smallskip
{\bf Bordism Lemma \Bor.} {\it
$\Omega_j(BO\left<5\right>)=\Z,\Z_2,\Z_2,\Z_{24},0,0,\Z_2,0$ according to
$j=0,1,2,3$, $4,5,6,7$.}

\smallskip
{\it Proof.}
We have $\Omega_j(BO\left<5\right>)\cong\Omega_j^{fr}$ for $j\le6$: for $j\le4$
because $BO\left<5\right>$ is 4-connected and for $j=5,6$ because
$\pi_j(BO\left<5\right>)\cong\pi_j(BO)\cong\pi_{j-1}(O)$.

Each map $Q\to BO\left<5\right>$ from a closed 7-manifold $Q$
bordant to a map from homotopy sphere $\Sigma$ [KM63].
(So $\Omega_7(BO\left<5\right>)=i_*\theta_7$, which is already sufficient for
the proof of the Primitivity Theorem \Pri.)
By [KM63, end of \S4] $\Sigma$ is a boundary of a parallelizable manifold,
and $BO\left<5\right>$-structure on $\Sigma$ extends to a
$BO\left<5\right>$-structure on this manifold.
Thus $\Omega_7(BO\left<5\right>)=0$.
\qed

\smallskip
{\it Proof of the injectivity of the attaching invariant
$a: E^7(S^4)\to\pi_4(G_3,SO_3)$ defined in \S3.}
Take embeddings $f,f':S^4\to\R^7$ such that $a(f)=a(f')$.
Then there exist framings $\xi$ and $\xi'$ of $\nu_f$ and $\nu_{f'}$ such that
$a(f,\xi)=a(f',\xi')$.
These framings define an orientation-preserving bundle isomorphism
$\varphi:\partial C_f\to\partial C_{f'}$.
Identify $\partial C_f$ and $\partial C_{f'}$ by $\varphi$.
Take any embedding $M_\varphi\to S^{16}$.

Let us set $C=S^2$ in the hypothesis of the Reduction Lemma \Red.
Using obstruction theory and the Fiber Lemma \fib we obtain a lifting
$l:M_\varphi\to BO\left<5\right>$.
Recall the definition of homotopy equivalences $h_f:C_f\to S^2$
and $h_{f'}:C_{f'}\to S^2$ from the construction of $a(f)$.
Since $a(f,\xi)\simeq a(f',\xi')$, we have $h_f\simeq h_{f'}$ on
$\partial C_f=\partial C_{f'}$.
Hence by the Borsuk Homotopy Extension Theorem $h_{f'}$ is homotopic to a map
$h':C_{f'}\to S^2$ such that $h'\varphi=h_f$ on $\partial C_f$.
Set $h=h_f\cup_\varphi h'$.

We have $\Omega_7(S^2\times BO\left<5\right>)=0$.
(This follows because in the Atiyah-Hirzebruch spectral sequence with
$E^2_{i,j}=H_i(S^2,\Omega_j(BO\left<5\right>))$ we have by the Bordism Lemma
\bor $E^2_{i,7-i}=0$.)
Hence $[h\times l]=0$.
Therefore by the Reduction Lemma \red $f$ is isotopic to $f'$.
\qed

\bigskip
{\bf Proof of the Primitivity Theorem \Pri.}

Denote $f'=f\#g$.
Since the normal bundle of $S^4$ in $S^7$ is trivial, there is an
orientation-preserving bundle isomorphism
$\varphi:\partial C_f\to\partial C_{f'}$ identical over $N_0$.
\footnote{There are several such $\varphi$ differing over $B^4$.
We shall prove that {\it each} of them works.}

Let us set $C=\C P^\infty=K(\Z,2)$ in the hypothesis of the Reduction Lemma
\Red.
Let $h_f$ be any map corresponding to $A_f[N]$ under the bijection
$H_5(C_f,\partial C_f)\to[C_f,\C P^\infty]$.
Define $h_{f'}:C_{f'}\to\C P^\infty$ analogously.
By the Complement Lemma \Com(b) $\pi_3(C_f)=\pi_3(C_{f'})=0$.
Hence $h_f$ and $h_{f'}$ are 4-connected.

%$$\minCDarrowwidth{5pt}\CD
%H_4(\partial C_f)@>> j\cong >  @>> exc\cong >
%H_4(\nu_f^{-1}N_0,\partial)\\
%@VV \varphi_* V @VV \varphi_* V @VV \varphi|_{N_0,*} V         \\
%H_4(\partial C_{f'})@>> j\cong > H_4(\partial C_{f'},\nu_{f'}^{-1}B^4)
%@>> exc\cong>  H_4(\nu_{f'}^{-1}N_0,\partial)\endCD.$$

Since $f'=f\#g$, we may assume that $fN_0=f'N_0$.
So by Remark \ru a weakly unlinked section for $f$ is a weakly unlinked section
for $f'$.
%Since $\varphi$ is identical over $N_0$, it follows that the $\varphi$-image
%of a weakly unlinked section for $f$ is a weakly unlinked section for $f'$.
Hence by the Section Lemma \Eu(a) (where $e$ and $j$ are isomorphisms)
$\partial A_{f'}[N]=\varphi_*\partial A_f[N]$.
The restrictions of $h_f$ and $h_{f'}\varphi$ to $\partial C_f$
correspond to $\partial A_f[N]$ and $\varphi^{-1}_*\partial A_{f'}[N]$
under the bijection $H_4(\partial C_f)\to[\partial C_f,\C P^\infty]$.
Hence these restrictions are homotopic.
Therefore by the Borsuk homotopy extension theorem $h_{f'}$ is homotopic to a
map $h':C_{f'}\to\C P^\infty$ such that $h'\varphi=h_f$ on $\partial C_f$.
Set $h:=h_f\cup_\varphi h'$.

Take any embedding $M_\varphi\to S^{16}$.
Since $C_f\subset S^7$, it follows that the (stable) normal bundle of $C_f$ is
trivial, so the classifying map $C_f\to BO$ of the normal bundle has a lifting
$C_f\to BO\left<5\right>$.
Obstructions to extending this lifting to $M_\varphi$ are in
$H^{i+1}(C_{f'},\partial C_{f'})\cong H_{6-i}(C_{f'})\cong H_{4-i}(N)$
with the coefficients $\pi_i(F)$.
Since $H_1(N)=0$, these obstructions are zeroes by the Fiber Lemma \Fib.
Thus there is a lifting $l_1:M_\varphi\to BO\left<5\right>$.

By the Reduction Lemma \red it remains to change $l_1$ to $l$ so that
$$[h\times l]=0\in\Omega:=
\Omega_7(\C P^\infty\times BO\left<5\right>)/i_*\theta_7.$$
Consider the Atiyah-Hirzebruch spectral sequence with
%\linebreak
$E^2_{i,j}=H_i(\C P^\infty,\Omega_j(BO\left<5\right>))$.
By the Bordism Lemma \bor among the groups $E^2_{i,7-i}$ the only non-trivial
ones are $E^2_{6,1}\cong\Z_2$ and $E^2_{4,3}\cong\Z_{24}$.
By [Te93, Proposition 1] the differential $E^2_{8,0}\to E^2_{6,1}$ is the
composition
$$\Z\cong H_8(\C P^\infty)\overset{\rho_2}\to\to H_8(\C P^\infty;\Z_2)
\overset{(\Sq^2)^*}\to\to H_6(\C P^\infty;\Z_2)\cong\Z_2.$$
Hence this differential is non-trivial.
Let $in:\C P^2\to \C P^\infty$ be the standard inclusion and
$p_5:BO\left<\infty\right>\to BO\left<5\right>$ the standard map.
Since $p_5^*:\Omega_3(BO\left<5\right>)\to\Omega_3^{fr}$ is an isomorphism,
the map
$$\Omega_3^{fr}\to\Omega\quad\text{defined by}\quad
(L:S^3\to BO\left<\infty\right>)\mapsto
(in\times p_5L:\C P^2\times S^3\to \C P^\infty\times BO\left<5\right>)$$
is an epimorphism.
Hence $[in\times p_5L]=-[h\times l_1]$ for some map
$L:S^3\to BO\left<\infty\right>$.

Since $h|_{C_f}$ is 4-connected, $(h|_{C_f})_*:H_4(C_f)\to H_4(\C P^\infty)$
is epimorphic.
Take $x\in H_4(C_f)$ such that $(h|_{C_f})_*(x)=[\C P^2]$.
Below we prove that

(*)\quad {\it each element of $H_4(C_f)$ is realized by an embedded
simply-connected 4-manifold.}

Thus $x$ is represented by a simply-connected
4-submanifold $X\subset C_f$.
Let $OX$ be a closed tubular neighborhood of $X$ in $M_\varphi$.
Denote by $D^3$ the fiber of the normal $D^3$-bundle $OX\to X$.

Define a map $l:D^3\to BO\left<5\right>$ so that $l\cup l_1|_{\overline D^3}$
would form a map homotopic to $p_5L$ (here $\overline D^3$ is $D^3$ with
reversed orientation).
Extend $l$ to a map $l:(M_\varphi-\Int OX)\cup D^3\to BO\left<5\right>$ as
$l_1$ outside $D^3$.
Obstructions to extending $l$ to $M_\varphi$ are in
$H^{i+1}(OX,D^3\cup\partial OX;\pi_i(F))$.
These obstructions are trivial by the Fiber Lemma \fib and
Lemma \clc below.
Hence $l$ extends to a lifting $l:M_\varphi\to BO\left<5\right>$.
Then $[h\times l]=[h\times l_1]+[in\times p_5L]=0$.
\qed

%\footnote{Note that $\Omega$ is isomorphic either to $\Omega_3^{fr}$ or to an index
%two subgroup of that group.}

%Among the differentials coming to or issuing out of the cell $(4,3)$
%On the line $i+j=7$ there could be only one
%another nontrivial differential $E^2_{4,3}\to E^2_{0,6}$.
%(This differential is in fact trivial because there is
%a splitting map
%$\Omega_6(BO\left<5\right>\times\C P^\infty)\to\Omega_6(BO\left<5\right>)$.)

\smallskip
{\bf Lemma \Clc.}
{\it If $OX\to X$ is a $D^3$-bundle over a closed 4-manifold $X$ such that
$H^1(X)=0$, then $H^{i+1}(OX,D^3\cup\partial OX)=0$ for $i=1,3$. }

\smallskip
{\it Proof.} For $i=1,3$ consider the exact sequence of triple
$(OX,D^3\cup\partial OX,\partial OX)$:
$$\minCDarrowwidth{5pt}\CD
H^i(OX,\partial OX)@>> > H^i(D^3\cup\partial OX,\partial OX)@>> >
 H^{i+1}(OX,D^3\cup\partial OX)@>> > H^{i+1}(OX,\partial OX)\\
@VV t V @VV \cong V @. @VV t V         \\
H^{i-3}(X) @>> r > H^i(S^3) @. @. H^{i-2}(X)
\endCD$$
Here $t$ are the Thom isomorphisms; we assume $H^a(X)=0$ for $a<0$.
For $i=1$ the statement is clear.
For $i=3$ the map $r$ is defined by restricting the normal bundle $OX\to X$
to a point, therefore $r$ is an isomorphism.
This, the commutativity and $H^1(X)=0$ imply that
$H^4(OX,D^3\cup\partial OX)=0$.
\qed

\smallskip
{\bf The Alexander Duality Theorem \Ale.}
{\it Let $f:N\to S^m$ be an embedding of a closed orientable $n$-manifold $N$.
The composition
$AD:H_{s+n-m+1}(N)\overset{\nu_f^!}\to\to H_s(\partial C_f)\overset i\to \to H_s(C_f)$
of the `preimage'
%homomorphism (from the Gysin exact sequence)
and the inclusion-induced homomorphisms is an isomorphism.}

\smallskip
This result
%version of Alexander duality
is apparently folklore, cf. [BRS76].
It holds because $AD$ is the composition of the preimage (=the Thom),
the excision and the boundary isomorphisms:
$$H_{s+n-m+1}(N)\overset{\widehat\nu_f^!}\to
\to H_{s+1}(S^m-\Int C_f,\partial C_f)\overset e\to
\to H_{s+1}(S^m,C_f)\overset\partial\to\to H_s(C_f).$$
Here $\widehat\nu_f:S^m-\Int C_f\to N$ is the $D^{m-n}$-normal bundle.
\footnote{In general neither $\nu_f^!$ nor $i$ is an isomorphism.
This `homology Alexander duality' is different from [Sk08', the Alexander Duality Lemma].
Isomorphism $AD$ coincides with the 'ordinary' Alexander duality, cf. [BRS76].
Indeed, $AD(x)$ spans the $(s+1)$-cycle $\widehat\nu_f^{-1}(x)$, so
the intersection (in $S^m$) of
$\widehat\nu_f^{-1}(x)$ with any $(n-(s+n-m+1))$-cycle in $fN$ is the same as
that of $x$ (in $N$), and the secondary linking coefficient (in $S^m$) of
$\widehat\nu_f^{-1}(x)$ with any torsion $(n-1-(s+n-m+1))$-cycle in
$fN$ is the same as that of $x$ (in $N$).}

\smallskip
{\it Proof of (*).}
By the Alexander Duality Theorem \ale each class in $H_4(C_f)$ can be
represented by $\nu_f^!y$ for some $y\in H_2(N)$.
The class $y$ is realizable by an embedded sphere with handles $M$
[Ki89, II, Theorem 1.1].
Hence $\nu_f^!y=[\nu_f^{-1}M]$.
We make $X:=\nu_f^{-1}M$ simply-connected by embedded surgery as follows.
If $\pi_1(X)\ne0$, then realize a generator of $\pi_1(X)$ by an embedding
$k:S^1\to X$.
Since $C_f$ is simply-connected, there is an extension $\overline k:D^2\to C_f$
of $k$.
By general position we may assume that $k\Int D^2\cap X=\emptyset$ and
$\overline k$ is an embedding.
A framing of $kS^1$ in $X$ can be extended to a triple of normal
linearly independent vector fields on $\overline kD^2$ because
$\pi_1(V_{5,3})=0$.
Thus $\overline k$ extends to an embedding $\hat k:D^2\times D^3\to C_f$ such
that $\hat k(\partial D^2\times D^3)\subset X$.
Let
$$X':=\left(X-\hat k(\partial D^2\times\Int D^3)\right)
\bigcup\limits_{\hat k(\partial D^2\times\partial D^3)}
\hat k(D^2\times\partial D^3),$$
$$\text{so that}\quad [X']=[X]\in H_4(C_f)\quad\text{and}\quad
\pi_1(X')\cong\pi_1(X)/\left<k\right>.$$
Continuing this procedure we get a simply-connected $X$.
\qed

%If $N$ is simply-connected, then we can make the above surgery inside $\partial C_f$, so $[X]=l[M]=[\nu_f^{-1}M]$.

\bigskip
{\bf Proof of the Almost Diffeomorphism Theorem \Adk.}

The `only if' part is easy (and is not used in the proof of the
Primitivity Theorem \Pri): take $W:=C_0\times I\sharp V$ (where $V$ is a
parallelizable manifold such that $\partial V=\Sigma$), take $B:=C_0\times BO$,
let $p$ be the projection onto the second factor, define $h:W\to C_0$ to be the
composition of the contraction of $V$, the natural diffeomorphism
$W/V\cong C_0\times I$ and the projection onto the first factor, and finally
let $\overline\nu:=h\times\nu$.

Let us prove the `if' part.
Since $C_0$ is simply-connected and $h|_{C_0}$ is 4-connected, $B$ is
simply-connected.
After surgery below the middle dimension we may assume that
$\overline\nu$ is 4-connected.
Since both $\overline\nu$ and $\overline\nu|_{C_k}$ are 4-connected, and
$W,B,C_k$ are simply-connected, the inclusion induces isomomorphisms
$H_i(C_k)\to H_i(W)$ and $H^i(W)\to H^i(C_k)$ for $i<4$.
Hence $H_5(W,C_{1-k})\cong H^3(W,C_k)=0$ by Poincar\'e-Lefschetz duality and
the exact sequence of pair $(W,C_k)$.
Consider the following diagram.
$$\minCDarrowwidth{5pt}\CD
@. @. K:=\ker\overline\nu_* @. @. \\
@. @. @VV \subset V @. @. \\
H_5(W,C_k)=0 @>> > H_4(C_k) @>> i_k > H_4(W) @>> j_k > V_k:=H_4(W,C_k) @>> > 0=H_3(C_k)\\
@. @. @VV \overline\nu_* V @. @.\\
@. @. H_4(B) @. @.
\endCD$$
Consider the following property ($P_{W,\overline\nu}$):
{\it there is a subgroup $U\subset K$ such that

$\bullet$ $U\cap U=0$,

$\bullet$ $j_k|_U$ is an isomorphism onto an additive direct summand
in $V_k$ for $k=0,1$, and

$\bullet$ the quotient $j_0U\times V_1/j_1U\to\Z$ of the intersection pairing
$\cap:V_0\times V_1\to\Z$ is unimodular. }
%$$V_1:=V_1\overset{j_1}\to\leftarrow H_4(W)\overset{j_0}\to\to V_0:=V_0.$$

\smallskip
{\it Completion of the proof of the `if' part of the Almost Diffeomorphism
Theorem \adk  under the assumption ($P_{W,\overline\nu}$).}
The form $\cap:K\times K\to\Z$ is even
%by the Divisibility Lemma \Di. \footnote{A simpler direct proof:
because
$$x\cap x=\left<w_4(W),x\right>=\left<\nu^*w_4,x\right>=\left<w_4,\nu_*x\right>
=\left<w_4,\pi_*\overline\nu_*x\right>=0\mod2,$$
where $\pi:BSpin\to BSO$ is the projection and $w_4\in H^4(BSO)$ is the
Stiefel-Whitney class.
So in [Kr99, p.725] we can take $\mu(x):=x\cap x/2$ for $x\in K$ (because 4 is even).
We have $Wh(\pi_1(B))=0$ and so an isomorphism is a simple isomorphism.
Hence the hypothesis on $U$ implies that
$\theta(W,\overline\nu)$ is `elementary omitting the bases' [Kr99, Definition
on p. 730 and the second remark on p. 732].
Thus the existence of the required diffeomorphism follows by the
$h$-cobordism theorem and [Kr99, Theorem 3 and the second remark on p. 732].
\qed

%The self intersection number of an embedded 4-sphere in an 8-manifold is mod 2
%equal to the evaluation of the Stiefel Whitney class $w_4$ on the homology
%class represented by the sphere (the self intersection number is the Euler
%class of the normal bundle, and this is mod 2 the Stiefel- whitney class of
%the normal bundle, which in turn (since $w_4$ is trivial for the tangent
%bundle of $S^4$) is the evaluation of $w_4$ on the homology class represented
%by the sphere.
%Bur our spherical classes all map to zero in $H_4(BSO)$, since they are in
%the kernel of the map to B. But $w_4$ is the pull back of the 4-Stiefel
%Whitney class in BSO, and by naturality the conclusion follows.
%Let us define 'elementary'. The pair $(W,\bar\mu)$ is {\it elementary} if

\smallskip
{\it Proof that the property ($P_{W',\overline\nu'}$) holds for some compact
8-submanifold $W'\subset S^{18}$ such that $\partial W'=M_\varphi$ and some
lifting $\overline\nu':W'\to B$ of the classifying map $W'\to BO$ of the
normal bundle.}
Since both $\overline\nu:W\to B$ and its restriction to $C_k$ are 4-connected,
maps $\overline\nu_*$ and $\overline\nu_*i_k$ on the diagram are surjective.
Hence by the Butterfly Lemma $j_k|_K$ is surjective.
Thus $j_k$ induces an isomorphism $j_k':V_k\to K/(K\cap\ker j_k)$.

If $j_0x=0$ then $x\cap y=j_0x\cap j_1y=0$  for each $y\in K$.
Since $j_1|_K$ is surjective and (by Poincar\'e-Lefschetz duality) the
intersection pairing $\cap:V_0\times V_1\to\Z$ is unimodular, the converse is
also true.
Hence
$$K\cap\ker j_0=\{x\in K\ |\ x\cap y=0\text{ for each }y\in K\}=
K\cap\ker j_1.$$
Therefore the bilinear form on $K':=K/(K\cap\ker j_0)$ defined by
$(a,b)\mapsto a\cap b=j_0a\cap j_1b$ is unimodular.
Since $\cap$ is even on $K$, the form on $K'$ is even.
Hence $\sigma(K')$ is divisible by 8.
Therefore there is the Kervaire-Milnor framed simply-connected 8-manifold
$V$ such that $\sigma(V)=-\sigma(K')$ and $\partial V$ is a homotopy sphere.
Stable framing on $V$ and map $\overline\nu$ give a map
$\overline\nu:W\natural V\to B$.
We have $\partial(W\natural V)=\partial W\#\partial V$.
Change of $W$ to the boundary connected sum $W\natural V$ has the effect of
making direct sum of the line
$V_0\overset{j_0}\to\leftarrow K\overset{j_1}\to\to V_1$ with
$A\overset{id}\to\leftarrow A\overset{id}\to\to A$, where $A$ is the
intersection form of $V$ so that $\sigma(A)=-\sigma(K')$.
Thus we may assume that $\sigma K'=0$.
Now standard argument (for new $K'$) implies that there is a submodule
$U'\subset K'$ satisfying the above three conditions with $K,U$ and $j_k$
replaced with $K',U'$ and $j_k'$.

(Indeed, there is $a_1\in K'$ such that $a_1\cap a_1=0$.
We can choose $a_1$ to be primitive.
By the unimodularity there is $b_1\in K'$ such that $a_1\cap b_1=1$.
The restriction of the intersection form to $\left<a_1,b_1\right>$ is
unimodular.
Hence $K/K'=\left<a_1,b_1\right>\oplus\left<a_1,b_1\right>^\perp$.
We may proceed by the induction to construct a basis
$a_1,\dots,a_s,b_1,\dots,b_s$.
Set $U':=\left<a_1,\dots,a_s\right>$.
Since elements of the basis are primitive, $U'$ is an additive direct summand.
It is also clear that $\cap:U'\times K'/U'\to\Z$ is unimodular.)

Since $H_i(C_1)\to H_i(W)$ are isomorphisms for $i<4$, we have $H_3(W,C_1)=0$.
Therefore by Poincar\'e-Lefschetz duality $V_0$ has trivial torsion.
Then $K'\cong V_0$ has trivial torsion.
Hence the projection $K\to K'$ has a right inverse $\psi$.
Then $U:=\psi U'$ satisfies the three conditions from the definition of property $(P_{W,\overline\nu})$.
\qed

%(The last sentence is true because when $\pi_1(B)=0$, by [Kr99, Definition on p. 725]
%and ??? we have $l_{4k}(\pi_1(B),w_1(B))\cong L_{4k}(0)\cong\Z$ and by
%[Kr99, Definition on p. 729] we have $\theta(W,\bar\mu)=\sigma(W)=0$.
%Therefore by [Kr99, Definition on p. 730] we obtain that $\theta(W,w')$ is elementary.)

\smallskip
Let us make some remarks
%(not used in the proof)
on higher-dimensional generalizations.

The proof shows that in the Almost Diffeomorphism Theorem \adk the dimension 4
can be replaced by $4q$ (and 3 by $4q-1$, \ 7 by $8q-1$, \ 8 by $8q$, \ 18 by $18q$).

%In some sense the 'converse' of the Reduction Lemma \red is true and can be
%used to construct invariants of embeddings and solve the following problems.

It would be interesting to describe $E^9(N)$ for closed connected 5-manifold $N$.
If $H_1(N)=0$, this reduces to description of the fibers of the Whitney invariant
$E^9(N)\to H_2(N;\Z_2)$ [MA] (which is surjective), i.e. of the orbits of the action
$\Z_2\cong E^9(S^5)\to E^9(N)$. Here we would need [Kr99, 5.ii].

%Madsen, Taylor and Williams, Theorem E mis-computes E^9(S^5)=0

It would be interesting to describe $E^8(N)$ for closed
connected smooth simply-connected 5-manifolds $N$.
The first step is to find $E^8(S^5)$
%which appears in exact sequences of Levine [L] and Haefliger [H]: a calculation is given in
[AC].
(Note that the calculations of $E^8(S^5)$ in previous versions of this paper is incorrect due to
a misreference to [Ha66]: $E^9(S^6)=0$ is not stated there and is presumably wrong.)

%See Levine's Corollary in 5.9 of the attached paper.  The 8-manifold M
% homotopy equivalent to S^6 \times S^2 obtained from sugery on a framed
%knot S^5 \times D^3 \to S^8.  If one is prepared to dig hard enough the
%relationship of this manifold to E^8(S^5) in explained in sections 3 and 5.

The Whitney invariant $W: E^{6k+4}(S^{2k+1}\times S^{2k+1})\to\Z_2\oplus\Z_2$
is injective because $ E^{6k+4}(S^{4k+2})=0$ [Mi72].
It would be interesting to know if the element $(1,1)$ is in its range (the
other elements of $\Z_2\oplus\Z_2$ are in the range [Sk06]).

\comment

One would need $E^8(S^5)\cong\Z_2$ (not stated in [Ha66, Mi72]).

\smallskip
{\it Proof of $E^8(S^5)\cong\pi_7(S^2)\cong\Z_2$.} From the exact sequence
$E^9(S^6)\to \theta_6^3\to \theta_6\to E^8(S^5)\to \theta_5^3\to\theta_5$
[Ha66, 4.20],
$E^9(S^6)=0$ [Ha66]
%this is not stated there and is wrong
and $\theta_5=\theta_6=0$, we obtain that
$\theta_6^3=0$ and $E^8(S^5)\cong\theta_5^3$.
Now the result follows by the following exact sequence (where
$\pi_n(G_3,SO_3)\cong\pi_{n+2}(S^2)$ analogously to the end of the proof of
Lemma \Gso.)
%[To62, tables]
$$\CD
\theta_6^3@>> >\pi_6(G_3,SO_3)@>> >P_6@>> >\theta_5^3@>> >\pi_5(G_3,SO_3)@>>
>P_5\\
@VV \cong V         @VV \cong V     @VV \cong V    @.            @VV \cong V    @VV \cong V\\
0         @.  \pi_8(S^2)\cong\Z_2  @.  \Z_2@.             @.
\pi_7(S^2)\cong\Z_2   @. 0
\endCD\quad\qed$$

% $C_6^4=C_6^3=0$ [Mi72], $C_5^4\cong\Z_2$ [Ha66] $C_5^3\cong\pi_7(S^2)\cong\Z_2$.
%This result can perhaps be recovered using the exact sequence [Ha66, 4.11] and not
%using $C_6^3=0$, if we prove that the stable suspension
%$\pi_6(\overline G_3)\to\pi_6(\overline G)$ is non-trivial.
%Note that
%$\pi_6(\overline G_3)\cong\pi_8(S^2)\cong\pi_6^S\cong\pi_6(\overline G)
%\cong\Z_2$
%but the suspension map $\pi_8(S^2)\to\pi_6^S$ is trivial.

\endcomment

%\newpage
\head 5. Appendix to \S3: an alternative proof of the Effectiveness Theorem \Eff.  \endhead

The idea is to construct directly the attaching invariant, which can perhaps be useful for generalizations.

\smallskip
{\bf Unlinked Framing Lemma.}
{\it If $N$ is simply-connected
%$\Sigma N\searrow\Sigma N_0$
and $f$ is PL compressible (see end of \S2), then

(a) there is a unique unlinked framing $\xi_f:N\times S^2\to\partial C_f$.

(b) each unlinked framing $N_0\times S^2\to\partial C_f$ extends to a unique
unlinked framing $N\times S^2\to\partial C_f$. }

\smallskip
The compressibility assumption is necessary for (a) but superfluous for (b).

The Unlinked Framing Lemma (a) is proved by first constructing an
unlinked framing $N_0\times S^2\to\partial C_f$ (the Extension Lemma \Exn(a)),
observing that it is unique (see the details below) and then using the Unlinked Framing Lemma (b).
%Proof of the latter is postponed and is an essential part of the proof of the Effectiveness Theorem.

\smallskip
{\it Proof of the uniqueness of a framing $N_0\times S^2\to\partial C_f$ for simply-connected $N$.}
Equivalence classes of framings $N_0\times S^2\to\partial C_f$ are in 1--1
correspondence with homotopy classes of maps $N_0\to SO_3$.
Obstructions to homotopy between maps $N\to SO_3$ are in
$H^i(N_0,\pi_i(SO_3))$.
By duality and since $\pi_2(SO_3)=0$ and $N_0/\partial N_0\cong N$, the latter
group is zero for each $i$.
% (for $i=1,3$ we use the condition $H_1(N)=0$, for $i=2$ we do not use it).
Therefore a framing $N_0\times S^2\to\partial C_f$ is unique.
\qed

\smallskip
{\it Proof of the Unlinked Framing Lemma (b).}
Since $f$ is compressible, $\sigma(N)=0$.
Therefore an unlinked framing $N_0\times S^2\to\partial C_f$ extends to a
framing $N\times S^2\to\partial C_f$ by the Extension Lemma \Exn(b).
Homotopy classes of such extensions $\xi$ and $\xi'$ differ by an element
$$d(\xi,\xi')\in H^4(N;\pi_4(SO_3))\cong\pi_4(SO_3)\cong\pi_4(S^2)\cong\Z_2.$$
Moreover, for fixed $\xi'$ the correspondence $\xi\mapsto d(\xi,\xi')$ is 1--1.
Identify the set $F$ of homotopy classes of such extensions with
$\pi_4(S^2)$ by this 1--1 correspondence.

Denote by $\xi_1$ the section formed by first vectors of $\xi$ and by
$l(\xi)$ the homotopy
class of the composition $N\overset{\xi_1}\to\to\partial C_f\subset C_f$.
Consider the following diagram:
$$\minCDarrowwidth{5pt}\CD @. @. \pi_4(S^2) @= F @. \\
@. @. @VV j_* V @VV l V @. \\
[\Sigma N,C_f]@>>r>[\Sigma N_0,C_f]@>> >\pi_4(C_f)@>>v>[N,C_f]@>> >[N_0,C_f]
\endCD$$
Here $j:S^2=S^2_f\to C_f$ is the inclusion and the bottom line is a segment
of the Barratt-Puppe exact sequence of $(N,N_0)$.
We do not know that the diagram is commutative but we know that
$l(\xi)-l(\xi')=vj_*d(\xi,\xi')$.
Hence $l(F)$ is a coset of $\im vj_*$.

Since $f$ is compressible, $\varkappa (f)=0$.
Hence by the Complement Lemma \Com(a) we have
$C_f\simeq C_{\varkappa (f)}\simeq S^2\vee(\vee_iS^4_i)$.
Identify these spaces.

Let us prove the existence of unlinked framing $N\times S^2\to\partial C_f$.
It suffices to prove that $l(\xi)\subset\im vj_*$ for a framing $\xi\in F$.
Since $\xi|_{N_0\times S^2}$ is unlinked, $l(\xi)|_{N_0}$ is homotopy trivial.
Therefore $vx=l(\xi)$ for some $x\in\pi_4(C_f)$.
We have $x=j_*y+z$ for some $y\in\pi_4(S^2)$ and $z\in\pi_4(\vee_iS^4_i)$.
Since $\xi|_{N_0\times S^2}$ is unlinked, $\varkappa (f)=0$ and $H_2(N)$ has no
torsion, $\xi_1|_{N_0}$ is weakly unlinked (in the sense of \S2).
So $(vx)_*=\xi_{1,*}:H_4(N)\to H_4(C_f)$ is trivial.
Thus $x_*:H_4(S^4)\to H_4(C_f)$ is trivial.
Therefore $z=0$.
So $l(\xi)=vx=vj_*y$.

Let us prove the uniqueness of unlinked framing.
 Like in the previous proof of the Effectiveness Theorem \Eff, $\Sigma N$ retracts to $\Sigma N_0$.
So $r$ is surjective.
Thus by exactness $v^{-1}(*)=0$.
Since $v$ extends to an action of the domain on the range, $v$ is injective.
Since $C_f\simeq S^2\vee(\vee_iS^4_i)$, we have that $j_*$ is injective.
So $vj_*$ is injective.
This implies the uniqueness.
\qed

 %We have that $\Sigma N$ retracts to $\Sigma N_0$ if $N$ is spin (that is, $w_2(N)=0$) and simply-connected [Mi58],
%or if $N=S^1\times S^3$, [Po85, Lecture 6], or if $N$ is a connected sum of manifolds with this property.
%It would be interesting to know whether $\Sigma N$ retracts to $\Sigma N_0$ for any spin 4-manifold $N$.

%This perhaps can be proved as follows.
%Represent $\Sigma N$ as a cell complex
%$\vee S^2_i\cup(\cup D^3_j)\cup(\cup D^4_k)\cup D^5$.
%We have
%$[D^5,D^4_k]=0$ because $H_4(N)\cong\Z$,
%$[D^5,D^3_j]=0$ because $w_2(N)=0$, so the intersection form of $N$ is even,
%hence $Sq^2:H^3(\Sigma N;\Z_2)\to H^5(\Sigma N;\Z_2)$ is trivial.
%$[D^5,S^2_i]=0$ because this is a suspension on some map $S^4\to S^1$.
%So $\Sigma N$ should retract to $\Sigma N-D^5=\Sigma N_0$.

\smallskip
{\it Proof of the Effectiveness Theorem \Eff.}
By the Unlinked Framing Lemma (a) there is a unique unlinked framing $\xi_f:N\times S^2\to\partial C_f$.
Take a retraction $r=r(\xi_f|_{N_0})$ given by the Retraction Lemma \Ret.
%Denote $A:=\xi(N_0\times S^2)$.
Since the framing $\xi$ is unlinked, the inclusion $\xi(N_0\times S^2)\to C_f$ extends to
a map $\xi(N_0\times S^2)\cup\Con(N_0\times*)\to C_f.$
By the Alexander duality and the Mayer-Vietoris sequence this map induces a homology isomorphism.
Hence by the relative Hurewicz Theorem this map is a homotopy equivalence.
Since projection $p:\xi(N_0\times S^2)\to S^2$ is null-homotopic on $N_0\times*$, projection $p$
extends to $\xi(N_0\times S^2)\cup\Con(N_0\times*)$.
Hence there is a retraction $r:C_f\to S^2_f$ whose restriction to $\xi(N_0\times S^2)$ is
the projection to $\xi(*\times S^2)=S^2_f$.

{\it The attaching invariant} $a_N(f)$ is the homotopy class of the composition
$$N\times S^2\overset{\xi_f}\to\cong\partial C_f\subset C_f\overset r\to \to S^2_f.$$
Clearly, $a_N(f)|_{N_0\times S^2}$ is homotopic to the projection onto the
second factor, $a_N(f)|_{x_0\times S^2}=\id S^2$ and $a_N(f)|_{N\times y_0}$
is null-homotopic.
Like in the previous proof $\Sigma N$ retracts to $\Sigma N_0$.
Therefore such maps are in canonical 1--1 correspondence with elements of $\pi_6(S^2)$
by the following Homotopy Lemma.
So we have $a_N(f)\in\pi_6(S^2)$.
For an embedding   $g:S^4\to S^7$ we have $a_N(f\#g)=a_{S^4}(g)+a_N(f)$.
This implies the Effectiveness Theorem \Eff.
\qed

\smallskip
{\bf Homotopy Lemma.}
{\it Let $n\ge3$ and $N$ be a closed $n$-manifold such that $\Sigma N$ retracts
to $\Sigma N_0$.
Denote by $X$ the set of homotopy classes of maps $a:N\times S^2\to S^2$
for which $a|_{N_0\times S^2}$ is homotopic to the projection onto the second
factor, $a|_{x_0\times S^2}=\id S^2$ and $a|_{N\times y}$ is null-homotopic.
% ($=y$ then $F_q$).
Then there is a canonical 1--1 correspondence $X\to\pi_{n+2}(S^2)$.}

\smallskip
{\it Proof.}
Maps $N\times S^2\to S^2$ for which $a|_{x_0\times S^2}=\id S^2$ can be
considered as maps $N\to G_3$.
Consider the Barratt-Puppe exact sequence of sets:
$$[\Sigma N;G_3]\overset r\to\to[\Sigma N_0;G_3]\to\pi_n(G_3)\overset v\to\to
[N;G_3]\to[N_0;G_3].$$
Since $\Sigma N$ retracts to $\Sigma N_0$, it follows that $r$ is surjective.
So by exactness $v^{-1}(*)=0$.
Since $v$ extends to an action of the domain on the range, $v$ is injective.
Therefore $v$ defines a 1--1 correspondence between $X$ and $\ker p$.
Here and below we use the notation from the diagram in the proof of
$\pi_4(G_3,SO_3)\cong\pi_6(S^2)$ in \S3.
Since $\partial$ factors through $\pi_n(SO_2)=0$, it follows that $\partial=0$.
Hence $\ker p=\im i\cong\pi_n(F_2)\cong\pi_{n+2}(S^2)$.
\qed

%\newpage

\enddocument